\documentclass[12pt]{amsart}
\usepackage{amssymb}
\numberwithin{equation}{section}
\let\cal\mathcal
\def\Ascr{{\cal A}}
\def\Bscr{{\cal B}}
\def\Cscr{{\cal C}}
\def\Dscr{{\cal D}}
\def\Escr{{\cal E}}
\def\Fscr{{\cal F}}
\def\Gscr{{\cal G}}
\def\Hscr{{\cal H}}

\def\Kscr{{\cal K}}
\def\Lscr{{\cal L}}

\def\Oscr{{\cal O}}

\def\Rscr{{\cal R}}

\def\Tscr{{\cal T}}

\let\blb\mathbb
\def\CC{{\blb C}}

\def \PP{{\blb P}}

\def \ZZ{{\blb Z}}

\def \NN{{\blb N}}

\def\opp{{\text{\upshape opp}}}

\def\Id{\operatorname{id}}
\def\pr{\mathop{\text{pr}}\nolimits}

\def\Lotimes{\overset{L}{\otimes}}

\def\Mod{\operatorname{Mod}}
\def\Vect{\operatorname{Vect}}
\def\mod{\operatorname{mod}}
\def\Gr{\operatorname{Gr}}

\def\QGr{\operatorname{QGr}}
\def\qgr{\operatorname{qgr}}

\def\gr{\operatorname{gr}}

\def\Qch{\mathop{\text{\upshape{Qch}}}}
\def\coh{\mathop{\text{\upshape{coh}}}}

\def\gr{\operatorname {gr}}
\def\Spec{\operatorname {Spec}}

\def\Ext{\operatorname {Ext}}
\def\Hom{\operatorname {Hom}}
\def\End{\operatorname {End}}
\def\RHom{\operatorname {RHom}}
\def\uRHom{\operatorname {R\mathcal{H}\mathit{om}}}

\def\im{\operatorname {im}}

\def\coker{\operatorname {coker}}
\def\ker{\operatorname {ker}}

\def\End{\operatorname {End}}

\def\add{\operatorname {add}}

\def\r{\rightarrow}

\DeclareMathOperator{\Tors}{Tors}
\DeclareMathOperator{\tors}{tors}

\def\dirlim{\mathop{\vtop{\baselineskip -100pt\lineskip 0.4ex
\setbox0\hbox{\upshape lim}\copy0\hbox to \wd0{\rightarrowfill}}}\limits}
\def\invlim{\mathop{\vtop{\baselineskip -100pt\lineskip 0.4ex
\setbox0\hbox{\upshape lim}\copy0\hbox to \wd0{\leftarrowfill}}}\limits}

\newtheorem{lemma}{Lemma}[section]
\newtheorem{proposition}[lemma]{Proposition}
\newtheorem{theorem}[lemma]{Theorem}
\newtheorem{corollary}[lemma]{Corollary}

\newtheorem{lemmas}{Lemma}[subsection]
\newtheorem{propositions}[lemmas]{Proposition}
\newtheorem{theorems}[lemmas]{Theorem}
\newtheorem{corollarys}[lemmas]{Corollary}

\theoremstyle{definition}

\newtheorem{definition}[lemma]{Definition}

\newtheorem{definitions}[lemmas]{Definition}

{

\newtheorem{step}{Step}

\newtheorem*{hypothesis}{Hypothesis}
\newtheorem*{properties}{Properties}
\newtheorem*{sublemma}{Sublemma}
}

\theoremstyle{remark}

\newtheorem{remarks}[lemmas]{Remark}

\newdimen\uboxsep \uboxsep=1ex
\def\uboxn#1{\vtop to 0pt{\hrule height 0pt depth 0pt\vskip\uboxsep
\hbox to 0pt{\hss #1\hss}\vss}}

\def\uboxs#1{\vbox to 0pt{\vss\hbox to 0pt{\hss #1\hss}
\vskip\uboxsep\hrule height 0pt depth 0pt}}

\usepackage{verbatim,amscd}

\def\epa{epaisse\ }
\def\Ab{\operatorname{Ab}}
\def\smd{\operatorname{smd}}
\def\add{\operatorname{add}}
\def\Ind{\operatorname{Ind}}
\title[Generators and representability of functors]{
Generators and representability of functors in commutative and noncommutative geometry}
\author{A. Bondal}
\address{%
Steklov Math. Institut, Gubkina 8, 117966, Moscow GSP-1,
Russia}
\email{bondal@mi.ras.ru}
\author{M. Van den Bergh}
 \address{Departement WNI,  Limburgs Universitair Centrum,
 3590 Diepenbeek, Belgium.}
  \email{vdbergh@luc.ac.be}
\thanks{The first author was partially supported by RFFI grants 97-01-00933
and 99-01-01144 and grant for support of leading scientific
groups N 0015-96085.}
\thanks{The second author is a senior researcher at the FWO. He also
  wishes to thank the Clay Mathematics Institute for material support
  during the preparation of this manuscript.}
\thanks{Most of the material in this paper was conceived during  visits of the
authors to the Max Planck Institute for Mathematics in Bonn and the
Mathematical Sciences Research Institute in Berkeley. The authors are
grateful to both
institutions for their kind hospitality.}
\thanks{The research described in this publication was made possible in part
by Award RM1-2089 of the US Civilian Researsh \& Development Foundation
for the Independent States of the Former Soviet Union (CRDF)}
\keywords{saturation, generators, representability, triangulated categories}
\subjclass{Primary 18E30}
\begin{document}
\begin{abstract}
We give a sufficient condition for an $\Ext$-finite triangulated category to be
saturated.
Saturatedness means that every contravariant cohomological functor of
finite type to vector spaces is representable. The condition consists in
existence of a strong generator. We prove that the bounded derived categories
of coherent sheaves on smooth proper commutative and noncommutative varieties
have strong generators, hence saturated. In contrast
the similar category for a smooth
compact analytic surface with no curves is \emph{not saturated}.

\end{abstract}
\maketitle
\section{Introduction and motivation}
\label{ref-1-0}
In this paper $k$ will be a field. Unless otherwise specified all
categories will be $k$-linear. If $\Dscr$ is a triangulated category
then a cohomological functor $H:\Dscr\r \Vect(k)$ is of finite type if
for all $A\in \Dscr$, $\sum_i \dim H(A[i])<\infty$.

This paper is inspired by the following result.
\begin{theorem}
\label{ref-1.1-1}
  Assume that $X$ is a regular projective variety over a field $k$. Let
  $\Dscr$ be the derived category of bounded coherent complexes on $X$. Then every
  contravariant cohomological functor of finite type on $\Dscr$ is
  representable.
\end{theorem}
This theorem was first announced in \cite{Bondal4}, but the proof
in {\em loc.cit}
%is rather complicated.
works only for functors which are homologically bounded with respect to the
standard $t$-structure. 
In the appendix to this paper we will give a
short proof of a generalization of Theorem \ref{ref-1.1-1} which
states that every contravariant cohomological functor of finite type
on the derived category of perfect complexes on a (possibly singular)
projective variety over a field is representable by a bounded complex
of coherent sheaves.

Theorem \ref{ref-1.1-1} is a motivation for the following definition \cite{Bondal4}.
\begin{definition} Assume that $\Dscr$ is $\Ext$-finite, i.e.\ $\sum_n
  \dim \Hom(A,B[n])<\infty$ for all $A,B\in \Dscr$. Then $\Dscr$ is
  \emph{(right) saturated} if every contravariant cohomological functor of
  finite type
  $H:\Dscr\r\Vect(k)$ is representable.
\end{definition}
Saturated triangulated categories are significant for non-commutative
algebraic geometry.  It can be argued that any definition of a
non-commutative proper scheme should give rise to an associated
saturated triangulated category which is the analogue of the bounded
derived category of coherent sheaves in the commutative case. Some
evidence for this point of view is given by \cite{Bondal3}.

One of the aims of this paper is to give an intrinsic criterion for
$\Dscr$ to be saturated. The central observation is that $\Dscr$
should be finitely generated in a suitable sense. If $E\in \Dscr$ then
we say that $E$ is a \emph{classical generator} for $\Dscr$ if $\Dscr$ is the smallest triangulated
subcategory of $\Dscr$ containing $E$ which is closed under summands.

If we define $\langle E\rangle_n$ to be the full subcategory of objects in
$\Dscr$ which can be obtained from $E$ by taking finite direct sums, summand, shifts
and \textbf{at most $n-1$ cones} then $E$ is a classical generator if
and only if $\langle E\rangle\overset{\text{def}}{=} \bigcup_n \langle
E\rangle_n=\Dscr$. We say that $E$ is a \emph{strong generator} for
$\Dscr$ if for some $n$ we have $\langle E\rangle_n=\Dscr$.

One of our main results is the following.
%We will extract the following corollary from Theorem \ref{ref-1.4-3}.
\begin{theorem}
\label{ref-1.5-4}
Assume that $\Dscr$ is $\Ext$-finite and has a strong
generator. Assume  in addition that $\Dscr$ is Karoubian (i.e.
every projector splits). Then
$\Dscr$ is saturated.
\end{theorem}
Let us give an idea of the proof of this theorem.
If $E$ is a classical generator then  using a method similar to the one used for the Brown
representability theorem one proves (see lemma \ref{ref-2.4.1-20})
%\begin{lemma}
%\label{ref-1.3-2}
that if $\Dscr$ is $\Ext$-finite and has a classical generator
  and $H: \Dscr\r
\Vect(k)$ is a contravariant cohomological functor of finite type then
there exists a directed system $(A_i)_{i\in \NN}\in \Dscr$ such that
$H=\dirlim \Hom(-,A_i)$.
%\end{lemma}
The final step in the traditional proof of the Brown
representability theorem consists in taking the homotopy
limit $\tilde{A}$ of the directed system  and proving that it
represents $H$.

Unfortunately in our setting $\tilde{A}$ is not defined,
because the definition of the homotopy limit depends on an infinite summation. To handle this problem,
we introduce  $n$-\emph{resolutions of} $H$ \emph{with respect
to a subcategory} $\Escr$ in $\Dscr$ (see \S\ref{ref-2.3-14}). Such a resolution is a directed system which gives a
good approximation for $H$ on the subcategory $\Escr$. At the price of
increasing $n$, it continues to be   a resolution with respect to $\Escr$ enlarged
by cones  and direct summands.
%So the main difficulty in proving saturatedness or related theorems is
%to control the $\Ind$-object $(A_i)_i$.  To handle this problem

%\begin{proof} Let $E\in \Dscr$ be a strong generator and let   $H: \Dscr\r
%\Vect(k)$ be a contravariant cohomological functor of finite
%type. Then according to Theorem \ref{ref-1.4-3}.2 $H$ will be a direct
%summand of $\Hom(-,Q_n)$ for some $n$. This direct summand corresponds
%to a projector on the functor $\Hom(-,Q_n)$. By Yoneda's lemma we
%obtain a corresponding projector on $Q_n$. By the assumption that
%$\Dscr$ is Karoubian, this projector corresponds to a summand of
%$Q_n$. It is easy to see that this summand represents $H$.
%\end{proof}

%Besides the above mentioned results we discuss several other side results
%which  may be interesting in their own right.
\medskip

We prove several results related to existence of generators and
(non)saturatedness for some types of categories of geometric and noncommutative
geometric origin.

\medskip

We discuss the existence of generators and strong generators for
  schemes.  In particular we prove that every quasi-compact, 
  quasi-separated scheme has a
  classical generator. In combination with a recent result of Keller
  \cite{Keller4} (see
  Theorem \ref{ref-4.1-31}) this shows that quasi-compact, quasi-separated schemes are affine in
  a DG- or  $A_\infty$-sense.
        We also prove that on a smooth scheme every  classical generator is a
strong generator. It follows that the bounded derived category of coherent
sheaves on a smooth proper scheme is saturated.

\medskip

We apply Theorem \ref{ref-1.5-4} to prove a
result which generalizes Theorem \ref{ref-1.1-1} to the non-commutative
case.  If $R$ is a (non-commutative) graded left coherent ring then there is a
natural category $\qgr(R)$ which is an analogue for coherent sheaves
on the projective scheme associated to a commutative graded ring. More precisely $\qgr(R)$ is the category of finitely
presented graded left $R$-modules modulo finite length modules.
In Theorem \ref{coherentring} we show that under appropriate
homological conditions on $R$ (which are analogous to those satisfied
by smooth projective varieties) the bounded
derived category of $\qgr(R)$ is saturated.
This application represents
our  original motivation for studying this subject.

\medskip

%We give a counter example to saturatedness in the analytic case.
In contrast with the case of algebraic varieties, we prove that the bounded
derived category of coherent sheaves (or, equivalently,  of complexes of sheaves with coherent
cohomology, see corollary \ref{6.3aa}) on a smooth compact analytic surface with no curves is not
saturated. The proof uses perverse coherent sheaves and a result from
\cite{ReVdB1}.

\medskip

Throughout the paper, if $\Escr$ is an abelian category then
$D^b(\Escr)$ and $D(\Escr)$ denote respectively the bounded and
unbounded derived category of $\Escr$. If $\Lambda$ is a ring or a
DG-algebra then $D(\Lambda)$ denotes its unbounded derived category.

The authors wish to thank Bernhard Keller, Maxim Kontsevich and Karen
Smith for useful conversations.
\section{Generators and resolutions in triangulated categories}
\label{ref-2-5}
\subsection{Generators}
\label{ref-2.1-6}
In this section we temporarily drop the assumption that our
triangulated categories are $k$-linear. In this section and the next
one we define various notions of generators for triangulated categories.

If $\Dscr$ is a triangulated category  then a triangulated subcategory $\Bscr$
of $\Dscr$ is called \epa (thick) if it is
closed under isomorphisms and direct summands. As was shown by Rickard
\cite{Ri} this is equivalent to  Verdier's original definition.

If $\Escr=(E_i)_{i\in I}$ is a set
of objects then  we say that $\Escr$ \emph{classically generates} $\Dscr$ if
the smallest \epa triangulated subcategory of $\Dscr$ containing
$\Escr$ (called the {\em epaisse envelope} of $\Escr$ in $\Dscr$) is equal to $\Dscr$ itself. We say that $\Dscr$ is
\emph{finitely generated} if it is classically generated by one object.

By the right orthogonal $\Escr^\perp$ in $\Dscr$ we denote the full
subcategory of $\Dscr$ whose objects $A$ have the property that
$\Hom(E_i[n],A)=0$ for all $i$ and all $n$. $\Escr^\perp$ is an \epa
subcategory of $\Dscr$. We say that $\Escr$  \emph{generates} $\Dscr$
if $\Escr^\perp=0$.
 Clearly if $\Escr$ classically generates $\Dscr$ then
it generates $\Dscr$, but the converse is false.

Assume now that $\Cscr$ is a triangulated category admitting arbitrary
direct sums. An object $B$ in $\Cscr$ is {\em compact} if $\Hom(B,-)$
commutes with direct sums. Let $\Cscr^c$ be the full subcategory of
$\Cscr$ consisting of compact objects. We say that $\Cscr$ is
\emph{compactly generated} if $\Cscr$ is generated by $\Cscr^c$.
The following is proved in \cite{Neeman}.
\begin{propositions} \label{ref-2.1.1-7} $\Cscr^c$ is Karoubian.
\end{propositions}
\begin{proof} (Sketch) Using the  standard limit argument one first proves
  that $\Cscr$ is Karoubian. Since a direct summand of a compact
  object is compact, this implies that $\Cscr^c$ is Karoubian.
\end{proof}
Then we have the following
result by Ravenel and Neeman \cite{Neeman3}.
\begin{theorems}
\label{ref-2.1.2-8}
Assume that $\Cscr$ is compactly generated. Then a set of objects
  $\Escr\subset \Cscr^c$ 
  classically generates
  $\Cscr^c$ if and only if it generates $\Cscr$.
\end{theorems}
\subsection{Strong generators}
\label{ref-2.2-9}
In what follows objects and subcategories will be considered in a fixed triangulated
category ${\cal D}$.

If ${\cal E}$ is a subcategory (or simply a set of objects), then we denote by ${\rm add}({\cal E})$ the minimal strictly full
subcategory in ${\cal D}$ which contains ${\cal E}$ and is closed under taking finite
direct sums and shifts. We denote by ${\rm smd}({\cal E})$ the minimal strictly full subcategory
which contains ${\cal E}$ and is closed under taking (possible) direct summands.

Following \cite{BBD}, one introduces a multiplication on the set of strictly full
subcategories. If ${\cal A}$ and ${\cal B}$ are two such subcategories, let ${\cal
A}\star {\cal B}$ be the  strictly full subcategory whose objects $X$
occur in a
triangle $A\to X\to B$ with $A\in {\cal A}$, $B\in {\cal B}$. This multiplication is
associative in view of the octahedron axiom. If ${\cal A}$ and ${\cal B}$ are closed
under direct sums and/or shifts, then so is ${\cal A}\star{\cal B}$.
\begin{lemmas}
\label{ref-2.2.1-10}
If ${\cal A}$ and ${\cal B}$ are closed under finite direct sums then:
\begin{itemize}
\item[i)]
$
{\rm smd}({\cal A})\star {\cal B}\ \subset \ {\rm smd}({\cal A}\star {\cal B}),
\ \ {\cal A}\star {\rm smd}({\cal B})\ \subset \ {\rm smd} ({\cal A}\star {\cal B});
$
\item[ii)]
$
{\rm smd}({\rm smd}({\cal A})\star {\cal B})={\rm smd}({\cal A}\star {\rm smd}({\cal
B}))={\rm smd}({\cal A}\star {\cal B}).
$
\end{itemize}
\end{lemmas}

\begin{proof} ii)  obviously follows from i). If $X \in {\rm smd}({\cal
A})\star {\cal B}$, then $X$  fits in a triangle $A_0\to
X\to B$ with $B\in {\cal B}$ and $A_0\oplus A_1=A$ for some  $A\in {\cal A}$.
If we add to this triangle the triangle $A_1\xrightarrow{\rm id}
A_1\to 0$, we get the
triangle $A\to X\oplus A_1 \to B$, which shows that $X\in {\rm smd} ({\cal
A}\star {\cal B})$. This proves  the first inclusion in i). The other
inclusion is similar.
\end{proof}
\begin{lemmas}
\label{epaisseenvelope}
The epaisse envelope of a strictly full triangulated subcategory
$\Ascr\subset \Bscr$ consists of summands of objects in $\Ascr$.
\end{lemmas}
\begin{proof}
By lemma \ref{ref-2.2.1-10} we have  \def\smd{\operatorname{smd}}
\[
\smd \Ascr\star \smd \Ascr\subset \smd (\Ascr\star\Ascr) =\smd \Ascr
\]
This proves the lemma.
\end{proof}

Now we define a new multiplication on the set of strictly full subcategories closed
under finite direct sums by the formula:
$$
{\cal A}\diamond {\cal B}\ = {\rm smd}({\cal A}\star {\cal B}).
$$
This multiplication is associative in view of  lemma \ref{ref-2.2.1-10}
and the associativity of $\star $:
$$
({\cal A}\diamond {\cal B})\diamond {\cal C}= {\rm smd}({\rm smd}({\cal A}\star {\cal
B})\star {\cal C})= {\rm smd}({\cal A}\star {\cal B}\star {\cal C})={\rm smd}({\cal
A}\star {\rm smd} ({\cal B}\star {\cal C}))={\cal A}\diamond ({\cal B}\diamond {\cal C}).
$$
Moreover the formula holds:
\begin{equation}
\label{ref-2.1-11}
{\cal A}_1\diamond {\cal A}_2\diamond \dots \diamond {\cal A}_n= {\rm smd} ({\cal
A}_1\star \dots \star {\cal A}_n).
\end{equation}

 Denote
$$
\langle \Escr\rangle_1= {\rm smd}({\rm add}(\Escr))
$$
$$
\langle \Escr\rangle_k=\langle \Escr\rangle_{k-1}\diamond \langle
\Escr\rangle_1={\rm
  smd}(\langle \Escr\rangle_1\star \dots \star \langle \Escr\rangle_1)\ (k \ {\rm factors}).
$$
$$
\langle \Escr\rangle =\bigcup_k\langle \Escr\rangle_k
$$
Thus $\langle \Escr \rangle $ is the epaisse envelope of $\Escr $ in $\Dscr$.
So $\Escr$ classically generates $\Dscr$ in the sense of \S\ref{ref-2.1-6} if and only if
$\langle \Escr\rangle=\Dscr$.
\begin{definitions}
We say that $\Escr$  {\em strongly generates}  ${\cal D}$ if ${\cal
  D}=\langle \Escr\rangle_k$, for some $k$.  We say that $\Dscr$ is
  \emph{strongly finitely generated} if it is strongly generated by
  one object.
\end{definitions}
In other words $\Escr$ strongly generates $\Dscr$ if we can get to any
object in $\Dscr$ from objects in  $\Escr$ by a universally bounded
number of cones.

Assume now that $\Cscr$ is a triangulated category admitting arbitrary
direct sums and let $\Escr$ be a set of objects in $\Cscr$.
We denote by $\overline{\rm add}({\cal E})$ the minimal strictly full
subcategory in ${\cal C}$ which contains ${\cal E}$ and is closed under taking arbitrary
direct sums and shifts. We define $\overline{\langle \Escr\rangle}_k$ in
the same way  as  $\langle \Escr\rangle_k$, but replacing $\add$ by
$\overline{\add}$.

\begin{comment}

Similarly we put $\overline{\langle
  \Escr\rangle}=\bigcup_k \overline{\langle\Escr\rangle}_k$. {}{}From the proof of the Brown representability theorem
\cite{Keller1,Neeman1} one easily obtains:

\begin{propositions}
\label{ref-2.2.3-12}
  (``Infinite devissage'') Assume that $\Escr$ consists of compact objects. Then
  $\Escr$  generates $\Cscr$ if and only if
  $\overline{\langle \Escr\rangle}=\Cscr$.
\end{propositions}
\begin{proof} Let $C\in\Cscr$. Then obviously $H=\Hom_\Cscr(-,C)$ is a
  representable functor. Now the proof of the Brown
  representability theoren shows that the representing object for $H$
  may be chosen inside $\overline{\langle\Escr\rangle}$. Since
  representing objects are unique it follows that $C\in
  \overline{\langle\Escr\rangle}$.
\end{proof}
\end{comment}
Analyzing the proof of Theorem \ref{ref-2.1.2-8} one obtains the
following statement:
\begin{propositions}
\label{ref-2.2.4-13}
 Assume that $\Escr$ consists of compact
  objects. Then $\overline{\langle\Escr\rangle}_k\cap
  \Cscr^c=\langle\Escr\rangle_k$.
\end{propositions}
\begin{proof}
The following is taken from Keller's  writeup of the proof of Theorem
\ref{ref-2.1.2-8} (see \cite[\S5.3]{Keller1}).
Let $M\in  \overline{\langle\Escr\rangle}_k\cap
\Cscr^c$.
Thus $M$ is
a summand of an object $Z\in
\overline{\langle\Escr\rangle}_{k-1}\star \overline{\add}(\Escr)$.

We now have a commutative diagram
\[
\begin{CD}
@. M @.\\
@. @VVV @.\\
Z_{k-1} @>>> Z @>>> Z' @>>>
\end{CD}
\]
where the lower row is a triangle
with $Z_{k-1}\in \overline{\langle\Escr\rangle}_{k-1}$, $Z'\in
\overline{\add}(\Escr)$.  Since $M$ is compact the composition $M\r
Z\r Z'$ factors through an object $M'$ in $\add(\Escr)$. {}From this we may
construct a morphism of triangles
\[
\begin{CD}
M_{k-1} @>>> M @>>> M' @>>>\\
@VVV @VVV @VVV \\
Z_{k-1} @>>> Z @>>> Z' @>>>
\end{CD}
\]
Repeating this construction we obtain a commutative diagram
\[
\begin{CD}
M_0 @>>> M_1 @>>> \cdots @>>> M_{k-1} @>>> M\\
@VVV      @VVV       @.         @VVV        @VVV\\
0  @>>> Z_1  @>>> \cdots @>>>  Z_{k-1} @>>> Z
\end{CD}
\]
By construction the cone of each of the upper maps lies in $\add(\Escr)$. Hence by the
octahedral axiom the cone $M''$
of the composition $M_0\xrightarrow{\alpha} M$ lies in $\add(\Escr)\star \cdots \star \add(\Escr)$
($k$ times).

Now consider the resulting commutative diagram
\[
\begin{CD}
M_0 @>\alpha>> M\\
@VVV @VVV\\
0 @>>> Z
\end{CD}
\]
The right vertical map is split and hence monic. It follows that
$\alpha$ is zero and hence $M$ is a summand of $M''$. This finishes
the proof.
\end{proof}

\subsection{Resolutions}
\label{ref-2.3-14}
As above $\Dscr$ is a triangulated category and $\Escr$  is a
subcategory of $\Dscr$. If $A\in \Dscr$ then we write $h_A$ for the
representable functor $\Hom(-,A)$.

Below we say that a directed system of abelian groups
$(G_i,d_i)_{i>0}$ is of order $n$ if the compositions of any $n$
consecutive transition maps is zero (following \cite{Kapranov1} we could
also say that $(G_i)_i$ is a complex of order $n$).

If $(F_i)_i$ and $(E_i)_i$ are of order $a$ and $b$ respectively and $ (F_i)_i\r (G_i)_i\r (E_i)_i$ is exact, then $(G_i)_i$ is easily seen to be of order $a+b$.
\begin{definitions}
\label{ref-2.3.1-15}
Assume that $H:\Dscr\r \Ab$ is a contravariant
  cohomological functor. Then an \emph{$n$-resolution} of $H$ with respect to
  $\Escr$ is a directed system of objects $(A_i)_{i>0}$ together
  with compatible natural transformations $\zeta_i:h_{A_i}\r H$ such
  that for any $E\in\Escr$, $p\in\ZZ$, $\zeta_i(E[p])$ is surjective and $\ker
  (\zeta_i(E[p]))_i$ is of order $n$. A \emph{resolution} of $H$ is a
  $1$-resolution.
\end{definitions}
\begin{lemmas}
\label{ref-2.3.2-16}
  If $(A_i)_i$ is an $n$-resolution of $H$ with
  respect to $\Escr$ then it is also an $n$-resolution with respect to
  $\langle \Escr \rangle_1$.
\end{lemmas}
The following key lemma is perhaps less obvious.
\begin{lemmas}
\label{ref-2.3.3-17}
  Assume that $(A_i)_i$ is an $a$-resolution of $H$ with
  respect to $\Escr\subset \Dscr$ and a $b$-resolution with
  respect to $\Fscr\subset\Dscr$. Then $(A_{i+b})_{i}$  is an $a+b$-resolution
 with respect to $\langle\Escr\rangle_1\diamond
  \langle\Fscr\rangle_1$.
\end{lemmas}
\begin{proof}
  We have $\langle\Escr\rangle_1\diamond \langle\Fscr\rangle_1=
  \smd(\smd\add\Escr\star\smd\add\Fscr)= \smd(\add\Escr\star\add\Fscr)$.
  In view of lemma \ref{ref-2.3.2-16} we may without loss of generality replace
   $\add(\Escr)$, $\add(\Fscr)$ by $\Escr,\Fscr$ and then again using
  lemma \ref{ref-2.3.2-16} it suffices to show that we have an $a+b$-resolution
 with respect to $\Escr\star\Fscr$.

Let $G\in \Escr\star\Fscr$. Then $G$ fits into a triangle $E\r G\r F$
with $E\in \Escr$, $F\in \Fscr$.
Define $K(U)_i$ and $C(U)_i$  as the directed
systems
given by the kernel and cokernel of $\Hom(U,A_i)\r H(U)$. We now look
at the following diagram:
{\tiny
\[
\strut\hskip -1cm
\begin{CD}
0 @. 0 @. 0 @. 0 @. 0\\
@VVV @VVV @VVV @VVV @VVV\\
K(E[p+1])_i @>>> K(F[p])_i @>>> K(G[p])_i @>>> K(E[p])_i @>>> K(F[p-1])_i \\
@VVV @VVV @VVV @VVV @VVV\\
\Hom(E[p+1],A_i)@>>> \Hom(F[p],A_i)@>>> \Hom(G[p],A_i)@>>>
\Hom(E[p],A_i)@>>>   \Hom(F[p-1],A_i)\\
@VVV @VVV @VVV @VVV @VVV\\
H(E[p+1])@>>> H(F[p])@>>> H(G[p])@>>>
H(E[p])@>>>   H(F[p-1])\\
@VVV @VVV @VVV @VVV @VVV\\
0 @>>> 0 @>>> C(G[p])_i @>>> 0 @>>>0 \\
@VVV @VVV @VVV @VVV @VVV\\
0 @. 0 @. 0 @. 0 @. 0
\end{CD}
\]
}
If we think of the spectral sequence associated to the acyclic double
complex formed by the two middle rows then we quickly
obtain the following:
\begin{enumerate}
\item $C(G[p])_i$ is a subquotient of $K(F[p-1])_i$.  It follows that the order
  of $C(G[p])_i$ is less than or equal to $b$. Since the
  transition maps in $C(G[p])_i$ are obviously surjective it follows that
  $C(G[p])_i=0$ for $i>b$.
\item There is an exact sequence:
\[
K(F[p])_i\r K(G[p])_i \r K(E[p])_i
\]
whence $K(G[p])_i$ has order $a+b$.\qed
\end{enumerate}
\def\qed{}\end{proof}
This lemma yields our main result.
\begin{propositions}
\label{ref-2.3.4-18}
 Assume that $(A_i)_i$ is a resolution of $H$  with
  respect to $\Escr\subset \Dscr$. Take $a\ge 1$. Then:
\begin{enumerate}
\item
  $(A_{ai})_i$ is a resolution of $H$ with respect to $\langle
  \Escr\rangle_a$.
\item $H$ is a direct summand of the representable functor $h_{A_{2a}}$
    when restricted to $\langle \Escr\rangle_a$.
\end{enumerate}
\end{propositions}
\begin{proof}
\begin{enumerate}
\item
{}{}From lemma \ref{ref-2.3.3-17} we obtain by induction that
  $(A_{i+a-1})_{i}$ is an $a$-resolution  with respect to $\langle
  \Escr\rangle_a$. The first element of this $a$-resolution is $A_a$. Hence
 $(A_{ai})_i$ is an
  honest resolution with respect to $\langle\Escr\rangle_a$.
\item
For $h_{A_a}(Z)$ we have an exact sequence:
\[
0\r \ker \zeta_a(Z)\r h_{A_a}(Z)\xrightarrow{\zeta_a} H(Z)\r 0
\]
The transition map $h_{A_a}(Z)\r h_{A_{2a}}(Z)$ kills $\ker \zeta_a(Z)$. Therefore
%By definition we have a map $\zeta_{2a}:h_{A_{2a}}\r H$. On the other
%hand if $Z\in \langle
%  \Escr\rangle_a$ then $\ker \zeta_{a+n}(Z)$ dies in
%  $h_{A_{2a+n}}(Z)$. Since $\zeta_{a+n}(Z)\r H(Z)$ is surjective by
%  the definition of resolution
we obtain
a map $\theta(Z):H(Z)\mapsto
  h_{A_{2a}}(Z)$, which is natural in $Z$.  It is easily seen that the
  composition $\zeta_{2a}(Z)\circ \theta(Z)$ is the identity on
  $H(Z)$. Therefore $H$ is a summand of $h_{A_{2a}}$ when
  restricted to $\langle \Escr\rangle_a$.\qed
\end{enumerate}
\def\qed{}\end{proof}
\subsection{Construction of resolutions}
\label{ref-2.4-19}
In this section $\Dscr$ is an $\Ext$-finite $k$-linear triangulated category.
\begin{lemmas}
\label{ref-2.4.1-20}
Let $E\in\Dscr$ and let
  $H:\Dscr\r \Vect(k)$ be a contravariant cohomological functor of
  finite type. Then $H$ has a resolution with respect to $E$.
\end{lemmas}
\begin{proof} This is proved in the same way as the Brown
  representability theorem \cite{Keller1,Neeman1}. For completeness let us repeat
  the construction of the resolution.

We start by taking $A_1=\oplus_{n} E[n]\otimes_k H(E[n])$. There is
an obvious canonical map $\zeta_1:h_{A_1}\r H$ which is surjective when
evaluated on $(E[n])_n$. Let $G=\ker \zeta_1$ and put
$B_1=\oplus_{n} E[n]\otimes_k G(E[n])$.  Then the composition
$h_{B_1}\r G\r h_{A_1}$ is by Yoneda's lemma given by a map
$\psi_1:B_1\r A_1$.
We now have a complex of
functors
\begin{equation}
\label{ref-2.2-21}
h_{B_1}\xrightarrow{h_{\psi_1}} h_{A_1}\xrightarrow{\zeta_1} H\r 0
\end{equation}
which is exact when evaluated on $(E[n])_n$.

Let $A_2$ be the cone of $B_1\xrightarrow{\psi_1} A_1$. Since $H$ is a
cohomological functor we have an exact sequence
\[
H(A_2) \r H(A_1)\r H(B_1)
\]
which by Yoneda's lemma translates into an exact sequence
\[
\Hom(h_{A_2},H)\r \Hom(h_{A_1},H)\r\Hom(h_{B_1},H)
\]
{}{}From \eqref{ref-2.2-21} it follows that  $\zeta_1$ is mapped to zero in
$\Hom(h_{B_1},H)$. Whence $\zeta_1$ lifts to a map $\zeta_2:h_{A_2}\r
H$. The fact that the composition
\[
h_{B_1}\r h_{A_1}\r h_{A_2}
\]
is zero combined with the exactness of \eqref{ref-2.2-21} on $(E[n])_n$ implies
that $\ker \zeta_1(E[n])$ is killed in $h_{A_2}(E[n])$. Thus it is
clear that if we repeat this construction we obtain a resolution
$(A_i,\zeta_i)_i$
of $H$ with respect to $E$.
\end{proof}

% To prove lemma \ref{ref-1.3-2} in the introduction it now
% suffices to note that  $H$ and $\dirlim h_{A_i}$ are cohomological
% functors which coincide on $(E[n])_n$ and hence on $\langle
% E\rangle$. So if $E$ is a classical generator then $H$ and $\dirlim h_{A_i}$
% coincide everywhere.

\begin{lemmas}
\label{ref-1.4-3}
Assume that $\Dscr$ is $\Ext$-finite.
  Let $H: \Dscr\r
\Vect(k)$ be a contravariant cohomological functor of finite type
 and let  $E$ be an arbitrary object in $\Dscr$. Then for all $n$ there exists an object $Q_n$ such
  that $H$ restricted to $\langle E\rangle_n$ is a direct summand of the representable functor
  $\Hom(-,Q_n)$.
\end{lemmas}
\begin{proof}
By lemma \ref{ref-2.4.1-20} $H$ as a resolution with respect to $E$. Then
in the notation of
Proposition \ref{ref-2.3.4-18} we may take $Q_n=A_{2n}$.
\end{proof}

\begin{proof}[Proof of Theorem \ref{ref-1.5-4}]
Let $E\in \Dscr$ be a strong generator and let   $H: \Dscr\r
\Vect(k)$ be a contravariant cohomological functor of finite
type. Then $\Dscr=\langle E\rangle_n$ for some $n$ and 
according to lemma \ref{ref-1.4-3} $H$ will be a direct
summand of $\Hom(-,Q_n)$. This direct summand corresponds
to a projector in the endomorphism ring of the functor $\Hom(-,Q_n)$. By Yoneda's lemma we
obtain a corresponding projector in $\End(Q_n)$. By the assumption that
$\Dscr$ is Karoubian, this projector corresponds to a summand of
$Q_n$. It is easy to see that this summand represents $H$.
\end{proof}
\subsection{A counter example}
In this section we show with a simple counter example that Theorem
\ref{ref-1.5-4}
is false if we only assume the existence of a generator (and not of a
strong generator).

Let $R=k[[x]]$ where $k$ is a field and let $\Escr$ be the category
of torsion $R$-modules. Let $S$ be the simple $R$-module. Then $S$ is
a generator for $\Dscr=D(\Escr)$. To see this, note that $\Escr$ is
hereditary and has enough injectives. So every object in $\Dscr$ is
the direct sum of its cohomology objects (see lemma
\ref{ref-5.2.8-49} below for a more general statement). Hence we
have to show that the right orthogonal of $S$ in $\Escr$ is
zero. Since $\Escr$ is closed under injective hulls in $\Mod(R)$ we
have  $\Ext^\ast_\Escr(S,M)=\Ext^\ast_R(S,M)$.
 If
$\Ext^\ast_R(S,M)$ is zero then $M$ is both $x$-torsion and uniquely
divisible by $x$. Hence $M=0$.

It is easy to see that the compact objects in $\Dscr$ are finite direct sums
of shifts of $S_n=R/x^n R$. {}{}From this it is clear that $S$ is not a
strong generator (the number of cones we need to reach $S_n$ depends
on $n$) and neither is any other object in $\Dscr^c$.

$\Dscr$ is also not saturated. Indeed if $E$ is the injective hull of
$S$ then $\Hom(-,E)$ defines a functor of finite type which is not
representable. This is a special case of the following more general
result proved in \cite{ReVdB1}.
\begin{lemmas}
\label{ref-2.1-22}
 Assume that $\Escr$ is an $\Ext$-finite
  abelian category of finite homological dimension in which every
  object has finite length.  Then  $D^b(\Escr)$  is
  saturated if and only if $\Escr\cong \mod(\Lambda)$ where $\Lambda$
  is a
  finite dimensional algebra of finite global dimension and
  $\mod(\Lambda)$ is the category of finite dimensional $\Lambda$-modules.
\end{lemmas}
In particular, the category $\Dscr$ we considered above cannot be saturated
since then it would have enough projectives, which is clearly not the case.

\section{Generators and strong generators for schemes.}
In this section we consider generators and strong generators for
certain types of schemes.
\subsection{Statement of results}
If $X$ is a scheme then by $\Qch(X)$ we will denote the category of
quasi-coherent $\Oscr_X$-modules.  If $X$ is noetherian then $\coh(X)$
is the category of coherent  $\Oscr_X$-modules. If $X$ is a ringed space
then $D(X)$ is the derived category of modules of $\Oscr_X$-modules
and if $X$ is a scheme then $D_{\Qch}(X)$ will be the derived category
of $\Oscr_X$-modules with quasi-coherent cohomology. It is clear that
$D(X)$ and $D_{\Qch}(X)$ admit arbitrary direct sums.

Quasi-coherent sheaves are well-behaved on \emph{quasi-compact,
  quasi-separated} schemes. Recall that a
quasi-compact scheme is a scheme that has a finite covering by affine
open subschemes and a quasi-separated scheme is a scheme such that the
intersection of any two affine open subschemes is
  quasi-compact. Actually it is sufficient to check this last
  condition on the affine opens of an arbitrary finite affine covering. 

A
  noetherian scheme is quasi-compact and quasi-separated.
If $X$ is quasi-compact quasi-separated then 
$\Qch(X)$ is a Grothendieck category \cite{thomasson}.

Our aim is to describe the category of compact objects in
$D_{\Qch}(X)$ for a quasi-compact, quasi-separated scheme.  Recall
that a complex on a scheme is said to be \emph{perfect} if it is
locally quasi-isomorphic to a bounded complex of vector bundles. In
particular a perfect complex is in $D_{\Qch}(X)$ and if $X$ is
quasi-compact then it is in $D^b_{\Qch}(X)$.

We will prove the following theorem.
\begin{theorems}
\label{ref-3.1-23}
Assume that $X$ is a quasi-compact, quasi-separated scheme. Then
\begin{enumerate}
\item The compact objects in $D_{\Qch}(X)$ are precisely the perfect
  complexes. 
\item $D_{\Qch}(X)$ is generated by a single perfect complex.
\end{enumerate}
\end{theorems}
Denote by $D_{\mathrm{perf}}(X)$ the category of perfect complexes on $X$.
\begin{corollarys}
\label{corperf}
If $X$ is a quasi-compact, quasi-separated 
then $D_{\mathrm{perf}}(X)$
is finitely generated.
\end{corollarys}
\begin{proof}
This follows from Theorem \ref{ref-3.1-23} and Theorem
\ref{ref-2.1.2-8}.
\end{proof}
We recall the following result for \emph{separated schemes}.
\begin{theorems} \label{eqicat}\cite{Lipman1,Neeman,Keller} If $X$ is quasi-compact and \emph{separated}
  then the canonical functor $D(\Qch(X))\r D_{\Qch}(X)$ is an
  equivalence.
\end{theorems}
This
result is false (even on the bounded derived categories) if we only
assume $X$ to be quasi-compact quasi-separated.  A counter example by
Verdier is
given in \cite[App I]{Illusie}.

If $X$ is smooth over a field (in particular separated) then using
Theorem \ref{eqicat} or directly it is easy to see that
$D^b(\coh(X))\cong D_{\mathrm{perf}}(X)$.
For smooth schemes we will prove the following result:
\begin{theorems}
\label{ref-3.2-24}
Assume that $X$ is smooth over a field (in particular separated). Then
$D^b(\coh(X))$ is strongly finitely generated. 
\end{theorems}
Presumably the last theorem is true under the weaker hypothesis that
$X$ is noetherian and regular.
\begin{corollarys}
\label{ref-3.3-25}
Assume that $X$ is smooth and proper over a field. Then $D^b(\coh(X))$
is saturated.
\end{corollarys}
\begin{proof}
This follows from Theorem \ref{ref-1.5-4} and Proposition
\ref{ref-2.1.1-7}.
\end{proof}
\begin{remarks} In characteristic zero one may give a different proof
  of Corollary \ref{ref-3.3-25} as follows. By Chow's lemma
  and Hironaka's theorem there is a birational dominant map $f:Y\r X$ such that
  $Y$ is projective and smooth. Since $X$ is smooth it has rational
  singularites and hence $Rf_\ast \Oscr_Y=\Oscr_X$.
        Then $f^*$ makes $D^b(\coh(X))$ into an admissible subcategory \cite{Bondal2}
  in $D^b(\coh(Y))$.
%It is then easy to see that
  In this situation, saturatedness of $Y$ (which follows from Theorem
  \ref{ref-1.1-1}) implies  saturatedness of $X$.

It is not clear to the authors if this proof can be generalized to
characteristic $p$.
\end{remarks}

Recently Bernhard Keller has proved the following result \cite{Keller4}
\begin{theorems}
\label{ref-4.1-31}
Let $\Escr$ be a Grothendieck category and
  assume that $\Ascr=D(\Escr)$ is generated by a compact
  object $E$. Then $\Ascr=D(\Lambda)$ where $\Lambda$ is a
  DG-algebra whose cohomology is given by $\Ext^\ast(E,E)$.
\end{theorems}

Combining this theorem with Theorem \ref{ref-3.1-23} we
find  the following corollary to our results
\begin{corollarys} \label{boundeddg} Assume that $X$ is a quasi-compact quasi-separated scheme. Then
  $D_{\Qch}(X)$ is equivalent to $D(\Lambda)$ for a suitable DG-algebra
  $\Lambda$ with bounded cohomology.
\end{corollarys}
\begin{proof} The fact that $\Lambda$ has bounded cohomology follows
  from lemma \ref{perfectbounded} below.
\end{proof}
Informally we may say that quasi-compact, quasi-separated schemes are
affine in a ``derived sense''.
\subsection{Extension of compact objects}
% The following reformulation of a theorem by Neeman models 
% the solution of the extension problem of perfect complexes on an open
% subscheme $j:U\r X$ given by Thomason in \cite{thomasson}.
% \begin{theorems} 
% \label{ravenelneeman}\cite[Thm 2.1]{Neeman3} Let $\Dscr$, $\Cscr$ be  
% triangulated categories admitting arbitrary direct sums with $\Dscr$
% compactly generated. Assume that $j^\ast:\Dscr\r \Cscr$,
% $j_\ast:\Cscr\r \Dscr$ is a pair of exact adjoint functors such that
% $j^\ast j_\ast$ is the identity and such that $j_\ast$ commutes with
% direct sums.  Then $\Cscr$ is compactly generated by $j^\ast\Dscr^c$. If
% in addition $\ker j^\ast$ is generated by compact objects in $\Dscr$
% then the full subcategory of $\Cscr^c$ given by the the objects in
% $j^\ast \Dscr^c$
% is triangulated  and its epaisse envelope is
% $\Cscr^c$ 
% \end{theorems}
% \begin{proof}
% \end{proof}

First recall the following.
\begin{theorems}  \label{ravenelneeman}\cite[Thm 2.1]{Neeman3} Let $\Dscr$ be compactly generated triangulated
  category admitting arbitrary direct sums and let $\Kscr$ be a
  triangulated subcategory which is closed under direct sums and which
  is in addition generated by objects which are compact in $\Dscr$.
  Put $\Cscr=\Dscr/\Kscr$. Then
\begin{enumerate}
\item $\Cscr$ admits arbitrary direct sums;
\item  $\Cscr$ is compactly generated;
\item $\Dscr^c$ maps to $\Cscr^c$ under the quotient functor;
\item   the induced functor $\Dscr^c/\Kscr^c\r \Cscr^c$ is fully
  faithful;
\item$\Cscr^c$ is the epaisse envelope of $\Dscr^c/\Kscr^c$.
\end{enumerate}
\end{theorems}
Assume that we are in the situation of the previous theorem and put $\Bscr=\Cscr^c$ and
let $\Ascr$ be the strict closure of $\Dscr^c/\Kscr^c$. Then $\Bscr$ is the
epaisse envelope of $\Ascr$. In this situation there is a
simple criterion to decide if an object  in $\Bscr$ lies
in $\Ascr$. This is contained in the following proposition.
\begin{propositions}
\label{object}
Let $\Ascr$ be a strictly full triangulated subcategory in a triangulated
category $\Bscr$ such that the epaisse envelope of $\Ascr$ is $\Bscr$.
Then an object $X$ in $\Bscr$ is in $\Ascr$
iff its representative $[X]\in K_0(\Bscr )$ belongs to the image
of $K_0(\Ascr )$.
\end{propositions}
We will give the proof below. In the situation of Theorem \ref{ravenelneeman}
this was proved in \cite{Neeman3}. In the case of schemes it is
\cite[Prop. 5.5.4]{thomasson}.

We immediately obtain the following corollary. 
\begin{corollarys}
\label{k0cor}
In the situation of Proposition \ref{object}  if $X\in \Bscr$ then
$X\oplus X[1]\in \Ascr$.
\end{corollarys}
The rest of this subsection is devoted to proving Proposition \ref{object}.

For an abelian  monoid $M$ with an operation $\oplus$,
denote by $F(M)$ the free abelian
group generated by elements of $M$ and
by $G(M)$ the quotient
of $F(M)$ by the subgroup $E(M)$ generated by elements
$[X\oplus Y]-[X]-[Y]$ taken for all pairs of elements $X, Y \in M$.

For an additive category $\Ascr $ denote by $G_+(\Ascr )$ the abelian
monoid with
elements the isomorphy classes of objects in
$\Ascr$ and with operation $\oplus $. We also use the notation $F(\Ascr )$,
$G(\Ascr )$, $E(\Ascr )$ for the corresponding groups $F(G_+(\Ascr ))$,
$G(G_+(\Ascr ))$, $E(G_+(\Ascr ))$.

The following lemma is classical and easy to prove. 
\begin{lemmas}\label{equing}
For two objects $X$ and $Y$ in an additive category $\Ascr$, $[X]=[Y]$ in
$G(\Ascr )$ iff there exists $Z\in \Ascr$ such that $X\oplus Z\cong Y\oplus Z$.
\end{lemmas}
If $\Ascr$ is a strictly full additive subcategory in $\Bscr$ then the natural morphism
$F(\Ascr )\to F(\Bscr )$ is obviously an embedding, which takes $E(\Ascr
)$ to $E(\Bscr )$. Thus, we may regard
$F(\Ascr )$, $E(\Ascr)$  as subgroups of $F(\Bscr )$, $E(\Bscr )$.

\begin{lemmas}\label{lf}
Let $\Ascr$ be a strictly full additive subcategory in an
additive category $\Bscr$
such that any object in $\Bscr$ is a direct summand of an object in $\Ascr$.
Then $E(\Bscr )\cap F(\Ascr )=E(\Ascr )$.
\end{lemmas}
\begin{proof}
Any element in $G(\Ascr )$ can be presented in the form $[X]-[Y]$ with $X,Y\in
\Ascr$.
Hence any element in $F(\Ascr )$ has the form $[X]-[Y]+v$ with $X,Y\in \Ascr$
and $v\in E(\Ascr )$.
Suppose this element is in $E(\Bscr )$. Since $E(\Ascr )\subset E(\Bscr )$, then
$[X]-[Y]\in E(\Bscr )$. Then by lemma \ref{equing}
there exists $Z\in \Bscr$ such that $X\oplus Z\cong Y\oplus Z$. By
 the assumption
we can find $Z'$ such that $Z\oplus Z'$ is in $\Ascr$. Then $X\oplus (Z\oplus
Z')\cong Y\oplus (Z\oplus Z')$. It follows that $[X]-[Y]\in E(\Ascr )$.
\end{proof}

The Grothendieck group ${K}_0(\Ascr )$
of a triangulated category $\Ascr$ is
the free abelian group generated by the isomorphy classes of
objects modulo the relations $[Y]=[X]+[Z]$ taken for all exact triangles
$X\to Y\to Z\to \cdots $. Denote by $I(\Ascr )$
the kernel of the natural homomorphism $G(\Ascr )\to {K}_0(\Ascr )$.

\begin{propositions}
\label{isomor}
Let $\Ascr$ be a strictly full triangulated subcategory in a
triangulated category $\Bscr$. Suppose that the epaisse envelope of $\Ascr$
coincides with $\Bscr$. Then:
\begin{itemize}
\item[(i)]
The induced homomorphism $G(\Ascr )\to G(\Bscr )$ is monic.
\item[(ii)]
The induced homomorphism $I(\Ascr )\to I(\Bscr )$ is an isomorphism.
\item[(iii)]
The induced homomorphism ${K}_0(\Ascr )\to {K}_0(\Bscr )$ is monic.
\end{itemize}
\end{propositions}
\begin{proof}
It is clear that $\Ascr $ satisfies the conditions of the last lemma.
A lifting to $F(\Ascr )$ of an element $x$ from the kernel of
$G(\Ascr )\to G(\Bscr )$
belongs to $E(\Bscr )$. Hence by the last lemma it is in
$E(\Ascr )$. Then $x$ is zero, and (i) is checked.

It follows from (i) that $I(\Ascr )\to I(\Bscr )$ is monic. Let us show it
is epic. The group $I(\Bscr )$ is the subgroup in $G(\Bscr )$ generated by
elements $[Y]-[X]-[Z]$ where $X\to Y\to Z$ is a triangle in $\Bscr$.
Find elements $X'$ and $Z'$ in $\Bscr $ such that $X'\oplus X$ and $Z\oplus
Z'$ are in $\Ascr $. Add the trivial triangles $X'\to X'\to 0$ and $0\to Z'\to
Z'$ to the primary triangle. Then we get the triangle:
$$
X'\oplus X\to X'\oplus Y\oplus Z'\to Z\oplus Z'.
$$
As the two extreme elements of the triangle are in $\Ascr $ then so is the middle
one. Hence $[X'\oplus Y\oplus Z']-[X'\oplus X]-[Z\oplus Z']$ is an element
in $I(\Ascr )$. Its image in $I(\Bscr )$ coincides with $[Y]-[X]-[Z]$ modulo
relations in $G(\Bscr )$. This proves (ii).

An element from the kernel of ${K}_0(\Ascr )\to {K}_0(\Bscr )$,
once lifted
to $G(\Ascr )\subset G(\Bscr )$, is in $I(\Bscr )$. Hence by (ii) it is in
$I(\Ascr )$. Then (iii) follows.
\end{proof}
\begin{proof}[Proof of Proposition \ref{object}]
In view of proposition \ref{isomor} we may regard $G(\Ascr )$ as a subgroup
of $G(\Bscr )$. Let $i: \Ascr \to \Bscr $ be the embedding functor. Denote
by ${\rm K_0}(i)$, $G(i)$ the corresponding homomorphisms of groups.
{}From the snake lemma and Proposition \ref{isomor}.(ii) it follows that the induced
homomorphism on cokernels ${\rm Coker}\ G(i)\to {\rm Coker}\ {\rm K}_0(i)$ is
monic. Hence the image in ${\rm K}_0(\Bscr )$ of an element $x$ in $G(\Bscr
)$ is in ${\rm K}_0(\Ascr )$ iff $x\in G(\Ascr )$.

%We need to show that $G_+(\Bscr )\cap G(\Ascr )=G_+(\Ascr)$.
Let us prove the following criterion for an object $X\in \Bscr $ to yield
an element in $G(\Ascr )$:
\begin{equation}\label{ga}
[X]\in G(\Ascr ) \Longleftrightarrow X\oplus A_1=A_2,
\end{equation}
for some $A_1, A_2$ in $\Ascr $.

Indeed, if $[X]\in G(\Ascr )$ then $[X]=[Y]-[Z]$ for some $Y,Z\in \Ascr $.
Therefore, $[X\oplus Z]=[Y]$. By lemma \ref{equing} there exists $W\in \Bscr
$, such that $X\oplus Z\oplus W\cong Y\oplus W$. By hypotheses we can
find $V\in \Bscr $, such that $U=W\oplus V\in \Ascr $. Then $X\oplus (Z\oplus
U)\cong Y\oplus U$. This proves (\ref{ga}).

But the right hand side of \eqref{ga} yields
a (split) exact triangle of the form $A_1\to A_2\to X$, i.e. $X\in \Ascr$.
\end{proof}

\subsection{Compact generators for  derived categories of
  quasi-coherent sheaves} 
%If $f:X\r Y$ is a morphism of ringed spaces
%then there is a right derived fuctor $Rf_\ast:D(X)\r D(Y)$ which was
%constructed by Spaltenstein \cite{Spaltenstein} using K-injective
%resolutions. 

Recall that an object in the homotopy category of
complexes is
K-injective if it is right orthogonal to the acyclic
complexes. Spaltenstein \cite{Spaltenstein} has proved that every complex of $\Oscr_X$-modules
on a ringed space $X$ has  a K-injective
resolution. Right derived functors are computed by evaluating the 
original functor on a K-injective resolution.

Most of the arguments below are based on Mayer-Vietoris type
 triangles. Let us indicate how these are constructed. 
 Assume
$X=U_1\cup U_2$ with $U_1,U_2$ open and put $U_{12}=U_1\cap U_2$. Let
$j_1$, $j_2$ and $j_{12}$ be the inclusions of $U_1,U_2$ and $U_{12}$
into $X$. By looking at stalks we see that we have a short exact
sequence in $\Mod(\Oscr_X)$:
\[
0\r j_{12!} \Oscr_{U_{12}}\r j_{1!} \Oscr_{U_1} \oplus j_{2!}
\Oscr_{U_2} \r \Oscr_X\r 0
\]
If $A\in D(X)$ then we obtain a triangle
\[
\uRHom(\Oscr_X,A)\r \uRHom(j_{1!} \Oscr_{U_1}, A) \oplus 
\uRHom(j_{2!} \Oscr_{U_2}, A)\r \uRHom(j_{12!} \Oscr_{U_{12}}, A)
  \r 
\]
For the definition of $\uRHom$ see \cite[Prop. 6.1]{Spaltenstein}

If  $A^\cdot$ is a K-injective complex on $X$ and $j:U\r X$ is an open
embedding 
then $j^\ast A^\cdot$ is K-injective on $U$. This follows from the existence
of the exact left adjoint $j_!$. {}From this we easily obtain
$\uRHom(j_!\Oscr_U,A)=Rj_\ast (j^\ast A)$. 
Hence we obtain a triangle 
\begin{equation}
\label{mv0}
A\r Rj_{1\ast}(j_1^\ast(A))\oplus Rj_{2\ast}(j_2^\ast(A))\r 
Rj_{12\ast}(j_{12}^\ast(A))\r
\end{equation}
{}From this triangle we may derive other Mayer-Vietoris type triangles
by applying suitable functors.
%Applying $R\Gamma$ 
%we find
%\begin{equation}
%\label{mv1}
%R\Gamma(X,A)\r R\Gamma(U_1,A)\oplus R\Gamma(U_2,A) \r
%R\Gamma(U_{12},A)\r
%\end{equation}
If $f$ is a map $X\r Y$ and  the
restrictions of $f$ to $U_1$, $U_2$, $U_{12}$ are denoted by $f_1$,
$f_2$, $f_{12}$ respectively then applying $Rf_\ast$ we obtain a triangle
\begin{equation}
\label{mv2}
Rf_\ast A\r Rf_{1\ast} (j_1^\ast(A))\oplus Rf_{2\ast} (j_2^\ast(A))\r 
Rf_{12\ast} (j_{12}^\ast(A))\r
\end{equation}
Let $E$ be another object in $D(X)$. Applying $\RHom(E,-)$ to
\eqref{mv0} we
find a triangle
\begin{equation}
\label{mv3}
\RHom(E,A)\r \RHom(j_1^\ast E, j_1^\ast A)\oplus \RHom(j_2^\ast E, j_2^\ast A)
\r \RHom(j_{12}^\ast E, j_{12}^\ast A))\r
\end{equation}
The Mayer-Vietoris triangles may be used in connection with the
following principle:
\begin{propositions} (Reduction principle)  Let $P$ be a property
  satisfied by some schemes. Assume in addition the following.
\begin{enumerate}
\item $P$ is true for affine schemes.
\item If $P$ holds for $U_1,U_2,U_{12}$ as above then it holds for
  $X$.
\end{enumerate}
Then $P$ holds for all quasi-compact quasi-separated schemes.
\end{propositions}
\begin{proof} (See the proof of \cite[Lemma 3.9.2.4]{Lipman})
  Let $X$ be quasi-compact, quasi-separated. Since $X$ has a finite
  affine cover it has a finite cover by quasi-compact separated
  schemes $X_1,\ldots, X_n$. We use induction on $n$. Put $U_1=X_1\cap
  \cdots \cap X_{n-1}$, $U_2=X_n$.  Being open subsets of $X$, $U_1$,
  $U_2$, $U_{12}$ are quasi-separated and by looking at affine covers
  of the $X_i$ we easily see that these subsets are also quasi-compact.
  Furthermore since $X_i\cap X_n$ is a subscheme of a separated
  scheme, it is itself separated. Hence $U_1$ and $U_{12}$ have
  coverings by $n-1$ quasi-compact separated schemes.
By induction we may assume now
  $n=1$, in other words, $X$ is separated. 

We now repeat the same argument with an affine cover $X=X_1\cup \cdots
\cup X_n$. Since $X$ is separated $X_i\cap X_n$ is affine and
induction on $n$ reduces us to the case $X$ affine and we are done.
\end{proof}
\begin{remarks}
It is easy to see that the class of quasi-compact, quasi-separated
schemes is the biggest class of schemes to which the reduction
principle is applicable (for all properties $P$).
\end{remarks}

A map $f:X\r Y$ between schemes is said to be
quasi-compact, resp.\ quasi-separated if for every affine open
$U\subset Y$ the inverse image of $U$ is quasi-compact, resp.\ 
quasi-separated. Quasi-compact and quasi-separated morphisms are
stable under composition and pullback.
\begin{theorems} \label{lipmanfundamental}
 \cite[Prop. 3.9.2]{Lipman} If $f:X\r Y$ is quasi-compact,
  quasi-separated then
\begin{enumerate}
\item $Rf_\ast$ maps $D_{\Qch}(X)$ into $D_{\Qch}(Y)$. 
\item If\/
  $Y$ is quasi-compact then the image of $D_{\Qch}(X)^{\le 0}$ lies in
  $D_{\Qch}(Y)^{\le N}$ for some $N$.
\end{enumerate}
\end{theorems} 
\begin{proof} Since this statement is crucial for what follows we
  sketch the proof.
  We may clearly
  assume that $Y$ is affine. Then by the Mayer-Vietoris triangle
  \eqref{mv2} and the reduction principle we may assume that $X$ is also
  affine.
  
  To prove (2) it is now clearly sufficient to prove that the image of
  $D_{\Qch}(X)^{\le 0}$ lies in $D_{\Qch}(Y)^{\le 0}$.  Let $A\in
  D_{\Qch}(X)^{\le 0}$.
    According to \cite[Prop 3.13]{Spaltenstein} $A$ has a so-called
  ``special'' K-injective resolution $I$. By construction $I$ is an inverse limit 
$\invlim_n I_n$ of left bounded injective resolutions of $\tau_{\ge
  -n} A$ such that $I_n\r I_{n-1}$ is split epi in every
  degree.

Now $f_\ast I$ is the ``sheaffication'' of $U\mapsto
\Gamma(f^{-1} (U), I)$ where $U$ runs through the affine opens of
$Y$. Note that $f^{-1}(U)$ is also affine. Hence it is sufficient to show
for all $V\subset X$ affine open that $\Gamma(V,I)=\invlim_n
\Gamma(V,I_n)$ is acyclic in degrees $> 0$. This is
clearly true for $\Gamma(V,I_n)$. Furthermore the map $\Gamma(V,I_n)\r 
\Gamma(V,I_{n-1})$ is surjective, and
a quasi-isomorphism in degree $\ge -n+1$. We can now conclude by
\cite[Lemma 0.11]{Spaltenstein} which guarantees under these
conditions that  $H^i(\invlim \Gamma(V,I_n))=\invlim
H^i(\Gamma(V,I_n))$ for all $i$.

Now we prove (1). 
Since we have an affine map it is clear that
$Rf_\ast$ maps $\Qch(X)$ to $\Qch(Y)$.  Hence to conclude it is
sufficient to prove that for $A\in D_{\Qch}(X)$ we have
  $H^i(Rf_\ast(A))=f_\ast(H^i(A))$. If $A\in D^+_{\Qch}(X)$ then this is
  clear by devissage. The case of arbitrary $A$ is handled by writing it as an
  extension
$\tau_{< -N}A\r A\r \tau_{\ge-N} A\r$ for $N\gg 0$.
\end{proof}
\begin{corollarys} 
Assume that $f:X\r Y$ is quasi-compact, quasi-separated. Then $Rf_\ast$
commutes with arbitrary direct sums on $D_{\Qch}(X)$.
\end{corollarys}
\begin{proof} This question is local on $Y$ so we may assume that $Y$
  is affine. Since a direct sums of injective resolutions is a complex
  of flabby sheaves, which are acyclic for $f_\ast$, and since in
  addition $f_\ast$ commutes with direct sums it is
  clear that $Rf_\ast$ commutes with arbitrary direct sums in
  $D(X)^{\ge -N}$ for all $N$.

Let $(A_i)_{i\in I}$ be a family of objects in $D_{\Qch}(X)$. Then according to Theorem \ref{lipmanfundamental}(2)  for $N$
large compared to $j$ 
we have the following sequence of equalities: $H^j(Rf_\ast(\oplus_i A_i))=H^j(Rf_\ast(\tau_{\ge -N}(\oplus_i
A_i)))=H^j(Rf_\ast(\oplus_i(\tau_{\ge -N} A_i)))=\oplus_i H^j(Rf_\ast(\tau_{\ge
  -N} A_i)) =\oplus_i H^j(Rf_\ast A_i)$. Thus we obtain that the canonical
map $\oplus_i H^j(Rf_\ast A_i)\r H^j(Rf_\ast(\oplus_i A_i))$ is a
quasi-isomorphism. 
\end{proof}
The following analogue of Serre's theorem is a special case of Theorem
\ref{eqicat}.  
\begin{corollarys} \cite{Lipman1,Neeman,Keller}
\label{affinecase}
Assume that $X=\Spec R$ is affine. Then the obvious
  functor $D(R)=D(\Qch(X))\r D_{\Qch}(X)$ has a quasi-inverse given by
  $R\Gamma(X,-)$.
\end{corollarys}
\begin{proof} It is easy to see that this amounts to showing that if
  $A\in D_{\Qch}(X)$ then $H^i(R\Gamma(X,A))=\Gamma(X,H^i(A))$. For $X$
  left bounded this is clear and we may reduce the general case to
  this using (the analogue for $R\Gamma$) of Theorem
  \ref{lipmanfundamental}(2) (in the 
  same way as in the previous corollary).
\end{proof}
Recall the following result \cite{Neeman}.
\begin{lemmas}
\label{perfectaffine}
If $R$ is a ring then the compact objects in $D(R)$ are precisely the perfect
complexes (bounded
complexes of finitely generated projective modules).
\end{lemmas}
\begin{lemmas} 
\label{perfectcompact}
If $X$ is quasi-compact quasi-separated and $E\in D_{\Qch}(X)$
  is perfect then $E$ is compact in $D_{\Qch}(X)$. 
\end{lemmas}
\begin{proof} In the notation of  \eqref{mv3}  it follows from the
  five-lemma that if $j_1^\ast E$, $j_2^\ast E$ and $j_{12}^\ast E$
  are compact then so is $E$.  By the reduction principle it is then
  sufficient to consider the affine case but this follows from lemma
  \ref{perfectaffine} and corollary \ref{affinecase}.
\end{proof}
The following lemma was needed for Corollary \ref{boundeddg}.
\begin{lemmas} 
\label{perfectbounded}
If $X$ is quasi-compact, quasi-separated, $E\in D_{\mathrm{perf}}(X)$ and
$F\in D^b_{\Qch}(X)$
 then
  $\RHom(E,F)$ is bounded.
\end{lemmas}
\begin{proof} This follows from \eqref{mv3} and the reduction principle.
\end{proof}

\begin{proof}[Proof of Theorem \ref{ref-3.1-23}] 
Our proof that $D_{\Qch}(X)$ is generated by a single perfect complex
 is a 
 modification of the proof of
  \cite[Prop. 2.5]{Neeman1}.   We
proceed by induction on the number of elements in an affine  covering
of $X$. The case
where $X$ itself is affine is obvious by Corollary \ref{affinecase}: the generating object is $\Oscr_X$. To perform the induction step
we consider the situation where $X$ has an open covering $U\cup Y$
with $Y$ quasi-compact and
$D_{\Qch}(Y)$  having a perfect generator $E$ and  $U=\Spec R$
being affine. Put $S=U\cap Y$ and let the inclusion maps be as in the
following diagram
\[
\begin{CD}
S @>\alpha>> U \\
@V \beta VV @VV \gamma V\\
Y @>> \delta > X
\end{CD}
\]
Let $V=X\setminus Y=U\setminus S$. Then $V$ is a closed subset of $U$
and $X$. Since $S$ is quasi-compact and $U$ is affine it follows that
$V$ is defined by a finite number of elements $f_1,\ldots, f_n\in R$.
Let $Q$ be the object in $D_{\Qch}(U)^c$ associated to the complex of
free $R$-modules $\otimes_i (R\xrightarrow{f_i} R)$.  
According to
\cite{Neeman} $Q$ is a compact generator for the kernel of the
restriction map $\alpha^\ast:D_{\Qch}(U)\r D_{\Qch}(S)$. 

Since the homology of $Q$ has support in $V$ it follows that
$R\gamma_\ast Q\mid Y=0$. Furthermore we have $R\gamma_\ast Q\mid U=Q$
(this holds for any $Q$ and any open immersion $U\subset X$). It
follows that $R\gamma_\ast Q$ is perfect. Furthermore from the
Mayer-Vietoris triangle \eqref{mv3} (with $U_1=Y$, $U_2=U$,
$E=R\gamma_\ast Q$ and $A=Z$) we obtain
\begin{equation}
\label{ref-3.2-28}
\Hom(R\gamma_\ast Q,Z)=\Hom(Q,Z\mid U)
\end{equation}
for any $Z\in D_{\Qch}(X)$.

Since $D_{\Qch}(U)$  is compactly generated 
and $\ker \alpha_\ast:D_{\Qch}(U)\r D_{\Qch}(S)$
is generated  by a compact object in $D_{\Qch}(U)$ it follows from
Theorem
\ref{ravenelneeman} and Corollary \ref{k0cor} that there exists $F\in D_{\Qch}(U)^c$ such that
$F\mid S=E'\mid S$ with $E'=E\oplus E[1]$. By Corollary \ref{affinecase} $F$
is a perfect complex.   The perfect complexes $F$ on $U$ and $E'$ on
$Y$ can be glued yielding a perfect complex on $X$ in the
following way. Define $P\in
D_{\Qch}(X)$ by the exact triangle
\[
P\r R\gamma_\ast F
\oplus R\delta_\ast E'
\r R{\delta\beta}_\ast(E'\mid S)\r
\]
(the middle arrow is the direct sum of the two obvious morphisms).
One can easily check that $\delta^\ast P=E'$, $\gamma^\ast P=F$ by applying $\delta^\ast$ and $\gamma^\ast$ to this triangle.
 Thus $P$ is perfect.

We claim that $C=P\oplus R\gamma_\ast Q$ is a
compact generator for $D_{\Qch}(X)$.

Assume that $Z$ is right orthogonal to $R\gamma_\ast Q$. Using
\eqref{ref-3.2-28}  we find that $Z\mid U$ is right
orthogonal to $Q$.
%Hence according to \cite{Neeman} we have
It follows $Z\mid
U\to R\alpha_\ast(Z\mid S)$ is an isomorphism (cf. \cite{Neeman})
and hence $R\gamma_\ast(Z\mid
U)=R(\delta\beta)_\ast(Z\mid S)$.  We then obtain from the
Mayer-Vietoris triangle \eqref{mv0}
that the map
\begin{equation}
\label{ref-3.3-29}
Z\r R\delta_\ast (Z\mid Y)
\end{equation}
 is an isomorphism.

Assume now in addition that $Z$ is right orthogonal to $P$. Then by
the isomorphism \eqref{ref-3.3-29} and adjointness we
obtain that $Z\mid Y$ is right orthogonal to $P\mid Y =E\oplus
E[1]$. Hence $Z\mid Y=0$. Again using the isomorphism \eqref{ref-3.3-29} we
obtain that $Z=0$. This finishes the proof of  the fact that 
$D_{\Qch}(X)$ is generated by a single perfect complex.

Now we will prove that all compact objects are perfect.  By lemma
\ref{perfectcompact} and Theorem \ref{ref-2.1.2-8} it follows that
every compact object is a direct summand of a perfect complex. But by
looking on an affine cover and invoking Corollary \ref{affinecase} we
see that a direct summand of a perfect complex is perfect.
\end{proof}

\subsection{Strong generators for smooth schemes}
In this section we prove Theorem \ref{ref-3.2-24}. The proof uses  an extension
of Beilinson's ``resolution of the diagonal'' argument. The idea for
this approach is due to Maxim Kontsevich.

%%%%%%%%%%%%%%%%% replaced segment %%%%%%%%%%%%%%%%%%%%%%%%%%%%%%%%%%%%%%%
%
%We need the following property of compact objects on schemes.
%\begin{propositions} \cite{Neeman1}
%\label{ref-3.2.1-30}
%Assume that $X$ is quasi-compact
%  quasi-separated. Then $E\in D(\Qch(X))$ is compact if and only if it
%  is locally perfect, i.e. if and only if it is locally quasi-isomorphic
%  to a finite complex of finitely generated free modules.
%\end{propositions}
%%%%%%%%%%%%%%%%%%%%%%%%%%%%%%%%%%%%%%%%%%%%%%%%%%%%%%%%%%%%%%%%%%%%%%%

\begin{lemmas} Let $f_1:X\r W$, $f_2:Y\r W$ be quasi-compact maps of
  quasi-compact schemes. Assume that $E$, $F$ are compact generators
  for $D_{\Qch}(X)$ and 
  $D_{\Qch}(Y)$. Then $E\boxtimes_W F$ is a compact
  generator for $D_{\Qch}(X\times_W Y)$.
\end{lemmas}
\begin{proof}
The fact that $E\boxtimes_W F$ is compact follows from Theorem
\ref{ref-3.1-23}. So we only need to show that $E\boxtimes_W
F$ is a generator. Assume that $Z$ is right orthogonal to $E\boxtimes_W
F$. Let $\pr_{1,2}$ be the projections of $X\times_W Y$ on the first and
the second factor. Since
\[
\Hom_{X\times_W Y}(E\boxtimes_W F, Z[m+n])=\Hom_{X\times_W Y}(L\pr^\ast_1
E,\uRHom_{X\times Y}(L\pr_2^\ast F,Z[m])[n])
\]
we deduce that $R\pr_{1\ast} \uRHom_{X\times_W Y}(L\pr_2^\ast F,Z[m])=0$ for $m$
arbitrary.

Now let $U$, $V$ and $T$ be  open affines in $X$, $Y$ and $W$ such
that $f_1(U)\subset T$, $f_2(V)\subset T$. We find
\[
0=\Gamma(U,R\pr_{1\ast} \uRHom_{X\times Y}(L\pr_2^\ast F,Z[m+n]))=\Hom_Y(F,R\pr_{2\ast}
(Z[m]\mid U\times_W Y)[n])
\]
{}From which we deduce $R\pr_{2\ast} (Z[m]\mid U\times_W Y)=0$.  Restricting
to $V$ yields $\Gamma(U\times_W V,Z[m])=\Gamma(U\times_T V,Z[m]) =0$.

Since $U,V,T,m$ are arbitrary and since $X\times_Z Y$ is covered by
the affine open sets $U\times_T V$  this implies $Z=0$ by Corollary \ref{affinecase}.
\end{proof}
\begin{proof}[Proof of Theorem \ref{ref-3.2-24}]
Assume that $X$ is smooth over the field $k$ and let $E$ be a compact
generator for $D_{\Qch}(X)$. Then $X\times X$ is smooth as
well and if $\Delta\subset X\times X$ is the diagonal then
$\Oscr_\Delta$ is compact by Theorem
\ref{ref-3.1-23}. Hence according to Theorem
\ref{ref-2.1.2-8} and the above lemma $\Oscr_\Delta\in\langle E\boxtimes E\rangle_k$ for
certain $k\in \NN$. Let $Z\in D_{\Qch}(X)$. Then $Z=R\pr_{
  1\ast}(\pr_2^\ast Z\Lotimes \Oscr_\Delta)$ and hence $Z\in \langle
R\pr_{1\ast} (\pr_2^\ast Z\otimes E\boxtimes E)\rangle_k=
\langle E\otimes R\Gamma(E\Lotimes Z)\rangle_k$. Since
$R\Gamma(E\Lotimes Z)$ is a complex of vector spaces we find that $Z\in
\overline{\langle E\rangle}_k$ and hence by Proposition
\ref{ref-2.2.4-13} $D_{\Qch}(X)^c=\langle E\rangle_k$. Since for smooth varieties we have $D_{\Qch}(X)^c=D^b(\coh(X))$. This
finishes the proof of Theorem \ref{ref-3.2-24}.
\end{proof}

\section{Derived categories for graded rings}
\label{ref-5-33}
In this section we will associate to a graded ring $R$ a category
$\QGr(R)$ which is a non-commutative analogue of the category of
quasi-coherent sheaves on a projective variety \cite{AZ, VO}. We will
prove that under appropriate homological conditions on $R$ the
category of compact objects $D(\QGr(R))^c$ in the derived category of
$\QGr(R)$ is strongly finitely generated and hence saturated.

If $R$ is coherent  then we may
also introduce a category $\qgr(R)$ which is analogous to the category of coherent
sheaves on a projective variety. Under the homological conditions alluded to above
we have  $D(\QGr(R))^c=D^b(\qgr(R))$.  Thus in this way we obtain a complete
non-commutative analogue to Theorem \ref{ref-1.1-1}.

\begin{comment}
\begin{itemize}
\item
Some stuff  for the case that $R$ is coherent.
\end{itemize}
\end{comment}
\subsection{Generalities}
\label{ref-5.1-34}
In this section we develop some rudiments of projective geometry for
graded rings. We begin with some of the standard material on functors
related to the category of graded $R$-modules. Since we do not assume initially that $R$ is noetherian or coherent we state some of the basic facts and give their proofs.

Below $R=k\oplus R_1\oplus R_2\oplus \cdots$ is a  connected
graded ring over a field $k$ with graded maximal ideal $m=\oplus_{n>0} R_n$.
Following \cite{VdB16} we assume throughout that
 $\dim \Ext^i(k,k)<\infty$ for all $i\ge 0$.
 In particular $R$ is
finitely presented. Note that this condition on $R$ is left-right symmetric.

$\Gr(R)$ denotes the category of graded left $R$-modules. For
$n\in\ZZ$, $\Gr(R)$ comes equipped with a shift functor $M\mapsto
M(n)$ where $M(n)$ is defined by $M(n)_j=M_{n+j}$.

 We will
write $\Ext^i_{\Gr(R)}(M,N)$ for the $\Ext$-groups in $\Gr(R)$ and
$\Ext^i_R(M,N)$ for the graded $\Ext$-groups $\oplus_n \Ext^i_{\Gr(R)}(M,N(n))$. Thus
$\Ext^i_{\Gr(R)}(M,N)=\Ext^i_R(M,N)_0$.

We say that $M\in\Gr(R)$ is torsion if it is locally finite
dimensional, or equivalently if for all $a\in M$ there exists $n$ such
that $m^na=0$.  Let $\Tors(R)$ denote the corresponding full
subcategory of $\Gr(R)$.  Since $R$ is finitely generated, $\Tors(R)$
is a localizing subcategory of $\Gr(R)$. Furthermore finitely
generated objects in $\Tors(R)$ are finite dimensional. Let
$\QGr(R)=\Gr(R)/\Tors(R)$. We define $\tau$ as the functor which
assigns to a graded $R$ module its maximal torsion module.  By
$\pi:\Gr(R)\r \QGr(R)$ we denote the quotient functor. By standard
localization theory $\pi$ is exact and commutes with colimits.  We
denote the (fully faithful) right adjoint to $\pi$ by $\omega$ and we
denote the composition $\omega\pi$ by $Q$. Since $\pi\omega$ is the
identity it follows $Q^2=Q$.

The shift functors $M\mapsto M(n)$ define
shift functors on $\QGr(R)$ for which we will use the same
notation. Finally we will write $\Oscr=\pi R$. Note that by
adjointness it follows that
\begin{equation}
\label{ref-5.1-35}
(R^i\omega M)_0=\Ext^i_{\QGr(R)}(\Oscr,M)
\end{equation}
for $M\in \QGr(R)$.
\begin{lemmas}
\label{ref-5.1.1-36}
For any
  directed system  $(N_i)_{i\in I}$ and for any $n$ we have
\[
\Ext^j_R(R/R_{\ge n},\injlim N_i)=\injlim_i \Ext^j_R(R/R_{\ge n},N_i)
\]
\end{lemmas}
\begin{proof}
  The fact that $\dim \Ext^i(k,k)<\infty$ implies that $R/R_{\ge n}$ has
  a graded resolution consisting of finitely generated free modules.
  {}{}From this fact the lemma is clear.
\end{proof}
\begin{lemmas}
\label{ref-5.1.2-37}
$R^i \tau$ commutes with filtered colimits (and hence with direct sums)
for all $i$.
\end{lemmas}
\begin{proof}
This follows from the description  \cite{stenstrom}
\begin{equation}
\label{ref-5.2-38}
R^i\tau =\dirlim_n \Ext^i_R(R/R_{\ge
  n},-)
\end{equation}
together with lemma \ref{ref-5.1.1-36}.
\end{proof}
\begin{lemmas}
\label{ref-5.1.3-39}
Assume that $T$ is torsion. Then
\[
R^i\tau T=0\qquad\text{for}\qquad i>0
\]
\end{lemmas}
\begin{proof}
By  lemma \ref{ref-5.1.2-37}  it suffices to prove
this
in the case that $T$
is finite dimensional. But then it is clear
from \eqref{ref-5.2-38}
 if we look at the degrees of the generators
of the modules occuring in a minimal free resolution of
$R/R_{\ge n}$.
\end{proof}
\begin{lemmas}
\label{ref-5.1.4-40}
$Q$ is given by
\[
QM=\dirlim_{n}\Hom_R(R_{\ge n},M)
\]
\end{lemmas}
\begin{proof} Standard localization theory \cite{stenstrom} tells us
\[
QM=\dirlim_{n}\Hom_R(R_{\ge n},M/\tau M)
\]
So we need to show $\dirlim_{n}\Ext^i_R(R_{\ge n},\tau M)=0$ for
$i\le 1$. The vanishing of $\dirlim_{n}\Hom_R(R_{\ge n},\tau M)$ is
obvious and since $\Ext^1_R(R_{\ge n},\tau M)=\Ext^2_R(R/R_{\ge n},\tau M)$ the
other  vanishing follows from lemmas \ref{ref-5.1.3-39} and \eqref{ref-5.2-38}.
\end{proof}
\begin{lemmas}
\label{longexact1}
For $M\in \Gr(R)$ there is a long exact sequence
\[
0\r \tau M  \r M \r Q M \r R^1\tau M \r 0
\]
and isomorphisms $R^iQ M =R^{i+1}\tau M$ for $i\ge 1$. In particular
$R^i Q$ vanishes on $\Tors(R)$ and commutes with filtered colimits.
\end{lemmas}
\begin{proof} This follows  from the long exact sequence obtained by
  applying $\dirlim{}_n \Hom_R(-,M)$ to the system of exact sequences $0\r
  R_{\le n} \r R \r R/R_{\ge n}\r 0$ and then invoking
lemma \ref{ref-5.1.4-40} and
  \eqref{ref-5.2-38}.
\end{proof}
\begin{lemmas}
\label{longexact2}
One has $R^i Q=R^i \omega \circ \pi$.
\end{lemmas}
\begin{proof} One has to show that if $E\in \Gr(R)$ is injective then
  $\pi E$ is acyclic for $\omega$. Let
\[
0\r \pi E\r F_0\r F_1\r F_2\r\cdots
\]
be an injective resolution of $\pi E$.  Since $\pi\omega$ is the
identity, applying $\omega$ to this sequence we get that
\begin{equation}
\label{ref-5.3-41}
0\r QE\r \omega F_0 \r \omega F_1\r \omega F_2\r \cdots
\end{equation}
is a complex with homology in $\Tors(R)$. 
Since $E\r QE$ has torsion kernel and cokernel, then applying $R^iQ$
to this morphism  we find (using the
vanishing of $R^iQ$ on torsion objects by lemma \ref{longexact1})
\[
R^i Q (QE)=R^i QE=\begin{cases}
QE&\text{if $i=0$}\\
0&\text{otherwise}
\end{cases}
\]
Similarly,  $R^jQ(\omega F_i)=0$ for $j>0$, given the fact that $\omega
F_i$ is injective by
adjointness.

Then the spectral sequence for hyper cohomology yields that
\eqref{ref-5.3-41}  becomes exact if we apply $Q$. 
Since $Q^2 E=QE$ and $Q\omega F_i=\omega F_i$ it follows that the
original sequence was already exact.
\end{proof}
\begin{lemmas}
\label{filtered}
$R^i\omega$ commutes with filtered colimits.
\end{lemmas}
\begin{proof} Let $(M_j)_j$ be a directed system in $\QGr(R)$. Then we have
\begin{multline*}
R^i\omega(\dirlim_j M_j)= R^i\omega(\dirlim_j \pi \omega M_j)
= (R^i\omega\circ \pi) (\dirlim_j \omega M_j)
=R^i Q(\dirlim_j \omega M_j)\\
= \dirlim_j R^i Q(\omega M_j)
=\dirlim_j (R^i \omega\circ\pi \circ \omega)(M_j)
=\dirlim_j R^i\omega(M_j)\qed
\end{multline*}
\def\qed{}\end{proof}
In the sequel we will make the following assumption on $\tau$:
\begin{hypothesis} $\tau$ has finite cohomological
  dimension, i.e. $R^n\tau =0$ for $n\gg 0$.
\end{hypothesis}
This hypothesis implies that $\omega$ and $Q$ also have finite
cohomological dimension by lemmas \ref{longexact1}, \ref{longexact2}.
 Furthermore using the methods in \cite{RD}
we may define the unbounded derived functors $R\tau$, $R\omega$, $RQ$
by means of acyclic resolutions. {}{}From the definition one easily deduces
the following properties:
\begin{properties}
\begin{enumerate}
\item
$R\tau$, $R\omega$, $\Ext^i(\Oscr,-)$, $RQ$, $R\pi=\pi$ commute with direct sums.
\item $R\omega$ is the right adjoint to $\pi$ and $\pi\circ R\omega=\Id$.
\item $R\tau$ is the left adjoint to the inclusion functor
  $D_{\Tors(R)}(\Gr(R))\r D(\Gr(R))$.
\item $R\omega\circ R\tau=0$.
\item $RQ=R\omega \circ \pi$.
\item
For $M\in D(\Gr(R))$ there is a triangle:
\begin{equation}
\label{ref-5.4-42}
R\tau M\r M\r R QM\r
\end{equation}
\end{enumerate}
\end{properties}
\subsection{Saturatedness}
\label{ref-5.2-43}
 In this section we will show that under suitable hypotheses the
 category $D(\QGr(R))^c$ is $\Ext$-finite and saturated.

We recycle the notations and assumptions of the previous
section. Recall that  if $M$ is a graded $R$-module
then Artin and Zhang \cite{AZ} say that $R$
satisfies $\chi(M)$ if $\dim \Ext^i_R(k,M)<\infty$ for all $i$. The
significance of this condition is the following
\begin{lemmas} \cite[Cor. 3.6(3)]{AZ}
\label{ref-5.2.2-44}
The following are equivalent.
\begin{enumerate}
\item $R$ satisfies $\chi(M)$.
\item For all $i$, $R^i\tau M$ is finite dimensional in every degree and in
  addition has right bounded grading.
\end{enumerate}
\end{lemmas}
Below we will say  that $R$ satisfies $\mu$ if it satisfies
$\chi(R)$.
\begin{lemmas}
\label{ref-5.2.3-45}
\begin{itemize}
\item[(i)]
$D(\QGr(R))$ is generated by $\{\Oscr(n)\}_{n\in\ZZ}$. 
\item[(ii)]
One has $D(\QGr(R))^c=\langle
  \Oscr(n)_{n\in\ZZ}\rangle$.
\end{itemize}
\end{lemmas}
\begin{proof}
Assume that $M\in D(\QGr(R))$ is right orthogonal to
$\{\Oscr(n)\}_{n\in\ZZ}$. Using adjointness this implies that $R\omega M$ is right
orthogonal to $R(n)$. Hence $R\omega M=0$, but then $0=\pi \circ
R\omega M=M$.

By property (1) above $\Oscr(n)$ is compact. Hence (ii)
follows from (i) together with Theorem \ref{ref-2.1.2-8}.
\end{proof}
\begin{corollarys} Assume that $R$ satisfies $\mu$. Then
$D(\QGr(R))^c$ is $\Ext$-finite.
\end{corollarys}
\begin{proof}
By  lemma \ref{ref-5.2.3-45} $D(\QGr(R))^c$ is classically generated by
$\{\Oscr(n)\}_{n\in\ZZ}$. Hence it suffices to prove that $\sum_i \dim
\Ext^i(\Oscr(m),\Oscr(n))$ is finite.  Now we have
$\Ext^i(\Oscr(m),\Oscr(n))=
\Ext^i(\Oscr,\Oscr(n-m))=(R^i\omega \Oscr)_{n-m}=(R^iQ R)_{n-m}$ by
\eqref{ref-5.1-35} and property (5). The corollary now
follows
from lemma \ref{ref-5.2.2-44} and
the triangle \eqref{ref-5.4-42}.
\end{proof}
%\begin{lemmas} \label{compact} One has $D(\QGr(R))^c=\langle
%  \Oscr(n)_{n\in\ZZ}\rangle$.
%\end{lemmas}
%\begin{proof}
%By property (1) above $\Oscr(n)$ is compact. Hence the current lemma
% follows from Lemma \ref{ref-5.2.3-45} together with Theorem \ref{ref-2.1.2-8}.
% \end{proof}
 We will now show that we can do better. In the rest of this section $\sigma_{\le }$, $\sigma_{\ge }$,  $\tau_{\le }$ and $\tau_{\ge }$ denote  respectively the ``stupid'' and ``smart'' truncations of complexes.
\begin{lemmas}
\label{ref-5.2.5-46}
Let $d$ be the cohomological dimension of
  $\omega$. Then there exists  a number $l\le 0$ such that $\Oscr(n)\in
  \langle \Oscr(k)_{l\le k\le 0}\rangle_{d+1} $ for all~$n>0$ (see \S\ref{ref-2.1.2-8} for notations).
\end{lemmas}
\begin{proof} Let
  $(F_{in})_{i\ge 0}$ be a minimal free resolution of $(R/R_{\ge
 n})(n)$ (where as usual $F_{in}$ is placed in complex degree $-i$). Clearly $F_{0n}=R(n)$ and the other $F_{in}$  are direct
 sums of $R(v)$'s with $v\le 0$.
 Put $Z_i=\ker(F_{in}\r F_{i-1,n})$. Then $\sigma_{\le -1}\sigma_{\ge -d-1}
  (\pi F_{\cdot n})$ represents an element of $\Ext^{d+1}(\Oscr(n),\pi
  Z_{d+1})$ which is zero by  \eqref{ref-5.1-35}.  Thus $\Oscr(n)$ is a direct
  summand of $\sigma_{\le -1}\sigma_{\ge -d-1}
  (\pi F_{\cdot n})$. This shows that  $\Oscr(n)\in
  \langle \Oscr(k)_{k\le 0}\rangle_{d+1} $.

To obtain the stronger conclusion of the proposition we have to bound above
the $u$ such that $R(-u)$ occurs in $\sigma_{\le -1}\sigma_{\ge -d-1}
   F_{\cdot n}$. That is we have to bound $u$ such that $k(u)$ occurs in
   $\Ext^i_R((R/R_{\ge
  n})(n),k)$ for $i\le d+1$. Since  $(R/R_{\ge
  n})(n)$ is an extension of $k(t)$, $0<t\le n$ we have to bound the
$k(u)$ occuring in  $\Ext^i_R(k(t),k)$ for $i\le d+1$ and $t> 0$. Since
$\Ext^i_R(k(t),k)=\Ext^i_R(k,k)(-t)$ such a bound is given by the maximal
$v$ such that $k(v)$ occurs in $\Ext^i_{R}(k,k)$ for $i\le d+1$.
\end{proof}

Now we discuss the case when $\QGr(R)$ has finite homological dimension.
Recall that  if $\Cscr$ is an abelian category then the homological
dimension of $\Cscr$ is the maximal $i$ such that there exist
$M,N\in\Cscr$ with the property that $\Ext^i_\Cscr(M,N)\neq 0$.
\begin{lemmas} Assume that $\QGr(R)$ has finite homological
  dimension. Then the functor $\tau$ has finite cohomological
  dimension.
\end{lemmas}
\begin{proof} This follows from combining Lemmas
  \ref{longexact1},\ref{longexact2} with
  \eqref{ref-5.1-35}.
\end{proof}

\begin{lemmas}
\label{ref-5.2.6-47}
Assume that $\QGr(R)$ has homological dimension $h<\infty$. Then for every
$M\in \QGr(R)$ one has $M\in \overline{\langle \Oscr(k)_k
  \rangle}_{h+1}$.
\end{lemmas}
\begin{proof}
This is proved by observing that if $M=\pi N$ then a sufficiently long
free resolution of $N$ splits in $\QGr(R)$. The same argument was used
in the proof of lemma \ref{ref-5.2.5-46}.
\end{proof}
\begin{lemmas}
\label{ref-5.2.7-48}
 Assume that $\QGr(R)$ has homological dimension
  $h$. Then one has $D(\QGr(R))=\overline{\langle  \Oscr(k)_k\rangle}_{2h}$.
\end{lemmas}
\begin{proof}
Let $U\in D(\QGr(R))$. It is easy to see that we can construct maps
$\alpha_i:Q_i\r U$ with the following properties:
\begin{enumerate}
\item $Q_i$ is a complex  consisting of (possibly infinite)
 direct sums of $\Oscr(k)$'s which
  starts in degree $ih+1$ and ends in degree $(i+1)h-1$.
\item $H^\ast(\alpha_i)$ is an isomorphism in homology in degrees
  $ih+2$ up to $(i+1)h-1$ and surjective in degree $ih+1$.
\end{enumerate}
Now put $Q=\oplus_i Q_i$, $\alpha=\oplus_i \alpha_i:Q\r U$  and let
$V$ be the cone of $\alpha$. We
find that $H^p(V)=0$ except when $h\mid p$. Invoking lemma \ref{ref-5.2.8-49}
below we find that $V=\oplus_i H^{ih}(V)$. By using lemma \ref{ref-5.2.6-47}
each of the $H^{ih}(V)$ can be produced by using at most $h$
cones. So the total number of cones we need is:
\[
h-2\text{(to produce $Q$)}+h\text{(to produce $V$)}+1\text{(to
  produce $U$ from $Q$, $V$)}=2h-1\qed
\]
\def\qed{}\end{proof}
The following lemma was used in the proof.
\begin{lemmas}
\label{ref-5.2.8-49}
Assume that $\Cscr$ is an abelian category which satisfies AB4 (exact
direct sums) and has enough injectives. Assume that the
homological dimension of $\Cscr$ is $h<\infty$ and let $V\in D(\Cscr)$
be a complex satisfying $H^p(V)=0$ unless      $h\mid p$. Then
$V=\oplus_i H^{ih}(V)$.
\end{lemmas}
\begin{proof}
  Write $H(V)=\oplus
H^{ih}(V)[-ih]$ (this sum exist since we have AB4 \cite{Neeman}).
We want to construct a quasi-isomorphism
$H(V)\r V$. To this end it is sufficient to construct maps $H^{ih}(V)[-ih]\r
V$ which induce isomorphisms on the $ih$'th cohomology. Since $\tau_{\le ih}
X\r X$ induces an isomorphism on $H^{ih}$, it is clearly sufficient to show
that the canonical map $\tau_{\le ih}V \r H^{ih}(V)[-ih]$ splits. {}{}From the
triangle
\[
\tau_{\le (i-1)h} V\r \tau_{\le ih} V \r H^{ih}(V)[-ih]\r
\]
we find that we have to show that
\begin{equation}
\label{ref-5.5-50}
\Hom(H^{ih}(V)[-ih],\tau_{\le (i-1)h}
V[1])=0
\end{equation}
 Now
 according to \cite[Thm 5.1, Cor. 5.3]{RD}, if $\Cscr$ has enough
 injectives and $\Hom(H^i(V),-)$ has finite cohomological dimension
 then we can compute $\Hom(H^i(V),-)$ (which is equal to 
$H^0(\RHom(H^i(V),-))$) by acyclic
 resolutions. It follows easily that an
 object in $D(\Cscr)^{\le -N}$  can be represented by an acyclic
 complex which is non-zero only in degree $\le -N+h$. This clearly
 implies \eqref{ref-5.5-50}.
\end{proof}

Some of the statements below will refer to the ring $R^\opp$. As a
rule we will decorate the corresponding notations by a
superscript ``$\opp$''.

\begin{lemmas}
\label{ref-5.2.9-51}
Assume that $\QGr(R)$ has  homological dimension $h<\infty$ and that $R$
satisfies $\mu^\opp$. Then for $n>0$, $\Oscr(-n)\in \langle
\Oscr(k) _{k\ge 0}\rangle_{h+1}$.
\end{lemmas}
\begin{proof} This is proved in a similar way as
  lemma \ref{ref-5.2.5-46}. We start with a minimal resolution of $(R/R_{\ge
  n})(n)^\opp$, dualizing we obtain a complex starting with $R(-n)$
  whose homology is finite dimensional (using the
  $\mu^\opp$-condition). Applying $\pi$ we obtain an exact sequence
  which start with $\Oscr(-n)$ and consists in higher degrees of
  direct sums of $\Oscr(k)$, $k\ge 0$. As in lemma \ref{ref-5.2.5-46}
  $\Oscr(-n)$ will be a direct summand of a truncation of length $h+1$
  of
  this exact sequence.
\end{proof}

\begin{lemmas}
\label{ref-5.2.10-52}
Assume that
%$\tau$
%and
$\tau^\opp$ has finite cohomological dimension
and that in addition
 $R$
  satisfies $\mu$ and $\mu^\opp$.  Assume furthermore that $\QGr(R)$
  has finite homological dimension.
Then there exist numbers $m\le 0$, $e\ge 1$  such that $\Oscr(n)\in
  \langle \Oscr(k)_{m\le k\le 0}\rangle_e $ for all~$n$.
\end{lemmas}
\begin{proof}
This follows by combining lemma \ref{ref-5.2.5-46} with lemma \ref{ref-5.2.9-51}.
\end{proof}

\begin{propositions} Assume that
%$\tau$ and
$\tau^\opp$ has finite
  cohomological dimension and that in addition $R$ satisfies $\mu$ and
  $\mu^\opp$. Assume furthermore that $\QGr(R)$ has finite
  homological dimension. Then the following holds.
\begin{enumerate}
\item $D(\QGr(R))=\overline{\langle \Oscr(k)_{a\le k\le 0}\rangle}_b$
  for some $a\le 0$, $b\ge 1$.
\item  $D(\QGr(R))^c=\langle \Oscr(k)_{a\le k\le 0}\rangle_b$.
\end{enumerate}
In particular the $\Ext$-finite triangulated category $D(\QGr(R))^c$
is strongly finitely generated.
\end{propositions}
\begin{proof}
(1) follows by combining lemma \ref{ref-5.2.10-52} with lemma
   \ref{ref-5.2.7-48}. (2) follows from  Proposition
   \ref{ref-2.2.4-13}.
\end{proof}
We can now finally prove the following theorem.
\begin{theorems}
\label{ref-5.2.12-53}
Under the hypotheses of the previous proposition
$D(\QGr(R))^c$ is saturated.
\end{theorems}
\begin{proof}
This follows Theorem \ref{ref-1.5-4}, Proposition
\ref{ref-2.1.1-7}
and the previous proposition.
\end{proof}
\def\uHom{\mathop{\Hscr\mathit{om}}}
\subsection{The case that $R$ is coherent}
Let $R$ satisfy the blanket assumptions made in the beginning of \S
\ref{ref-5.1-34} and assume that $R$ is left graded coherent. In other words the
kernel of a graded map between two free graded $R$ modules of finite rank is
finitely generated. Let $\gr(R)$ be the category of finitely presented
graded $R$-modules. Since $R$ is coherent this is an abelian
category.
%\begin{lemmas} Assume that $\tau$ has finite cohomological
%  dimension. Then $R$ satisfies $\mu$ if and only if it satisfies
%  $\chi(M)$ for all finitely generated $R$-modules.
%\end{lemmas}
%\begin{proof}
%  This is clear if we replace $M$ by a minimal free resolution and
%  then invoke lemma \ref{XXX}.
%\end{proof}
%Following \cite{AZ} we say that $R$ satisfies $\chi$ if it satisfies
%$\chi(M)$ for all finitely generated $R$-modules.

Put $\tors(R)=\gr(R)\cap \Tors(R)$. Then $\tors(R)$ consists
of the finite dimensional graded $R$-modules. We put
$\qgr(R)=\gr(R)/\tors(R)$. It is easy to see that the obvious functor $\qgr(R)\r \QGr(R)$ is
fully faithful.
\begin{lemmas}
\label{commutefiltered}
Let $M\in \qgr(R)$. Then $\Ext^i_{\qgr(R)}(M,-)$ commutes with filtered colimits.
\end{lemmas}
\begin{proof} By Lemma \ref{filtered} and \eqref{ref-5.1-35}  this is clearly true if
  $M=\Oscr(n)$ and it is a tautology if $i<0$. To treat the general we construct a short exact
  sequence
\[
0\r N\r F\r M\r 0
\]
where $F$ is a finite sum of shifts of $\Oscr(n)$.
Let $(T_j)_j$ be a directed system. We now have the following
commutative diagram
{\tiny
\[
\begin{CD}
\dirlim_j\Ext^{i-1}(F,T_j) @>>> \dirlim_j\Ext^{i-1}(N,T_j) @>>>\dirlim_j \Ext^i(M,T_j) @>>>
\dirlim_j\Ext^i(F,T_j)@>>> \dirlim_j\Ext^i(N,T_j)\\
@V\alpha VV @V\beta VV @V\gamma VV @V\delta VV @V\epsilon VV\\
\Ext^{i-1}(F,\dirlim_jT_j) @>>> \Ext^{i-1}(N,\dirlim_jT_j) @>>> \Ext^i(M,\dirlim_jT_j) @>>>
\Ext^i(F,\dirlim_jT_j)@>>> \Ext^i(N,\dirlim_jT_j)
\end{CD}
\]
}
$\alpha$ and $\delta$ are isomorphisms by the above
discussion. Furthermore we may assume by induction that $\beta$ is an
isomorphism. It now follows by diagram chasing that $\gamma$ is
monic. Then, replacing $M$ by $N$ we find that $\epsilon$ is also
monic. Performing another diagram chase yields that $\gamma$ is also epic.
\end{proof}
\begin{lemmas} Assume $\QGr(R)$ has finite cohomological
  dimension. Then $D(\QGr(R))^c=D^b_{\qgr(R)}(\QGr(R))$.
\end{lemmas}
\begin{proof} By Lemma \ref{ref-5.2.3-45}  $D(\QGr(R))^c$ is
  classically generated by $\{\Oscr(n)\}_{n\in\ZZ}$. Since
  $\Oscr(n)\in \qgr(R)$ this proves one inclusion.

To prove the other inclusion we have to show that every $M\in \qgr(R)$
is compact. This follows easily from the fact that by hypotheses
$\Ext^i(M,-)$ has finite cohomological dimension combined with  Lemma \ref{commutefiltered}.
\end{proof}
To conclude we give an alternative description of
$D^b_{\qgr(R)}(\QGr(R))$.
\begin{lemmas} The canonical functor $D^b(\qgr(R))\r
  D^b_{\qgr(R)}(\QGr(R))$ is an equivalence.
\end{lemmas}
\begin{proof}
 According to the dual version of \cite[1.7.11]{KSI} it is sufficient to
  prove the following result: if $B\r C$ is an epimorphism in $\QGr(R)$
  with $C\in \qgr(R)$ then there exists a map $D\r B$ with
  $D\in \qgr(R)$ such that the composition $D\r B\r C$ is an
  epimorphism.

The map $B\r C$ is obtained from a map $\theta:B_0\r C_0$ in $\Gr(R)$ with
$C_0\in \gr(R)$. But then the cokernel of $\theta$ is finite
dimensional and hence without loss of generality we may assume that
$\theta$ is epic. Since $C_0$ is finitely generated we may select a
finitely generated graded submodule $D_0$ of $B_0$  which contains inverse
images of the generators of $C_0$. This proves what we want.
\end{proof}
Combining everything we now obtain:
\begin{theorems}
\label{coherentring}
Let $R$ be a graded left coherent ring which satisfies the following
  hypotheses.
\begin{enumerate}
\item $\dim \Ext^i(k,k)$ is finite dimensional for all $i$.
\item $R$ satisfies $\mu$ and $\mu^{\opp}$.
\item $\tau^{\opp}$ has finite cohomological dimension.
\item $\QGr(R)$ has finite homological dimension.
\end{enumerate}
Then $D^b(\qgr(R))$ is $\Ext$-finite and saturated.
\end{theorems}

\section{Derived categories of analytic surfaces}
We have shown in Corollary \ref{ref-3.3-25} that if $X$ is a
smooth proper algebraic variety over a field $k$ then $D^b(\coh(X))$ is saturated.
Since smooth proper algebraic varieties and compact analytic manifolds have
 similar properties it is a natural question to ask if this result
remains true if we assume that $X$ is compact analytic. In this
section we show that the answer to this question is negative.

\subsection{Serre functors}
Let $X$ be a connected compact complex analytic manifold of
dimension $n$.  Write $D^b_{\mathrm{coh}}(X)$ for the bounded
derived category of sheaves of $\Oscr_X$-modules with coherent
cohomology. We first prove that $D^b_{\mathrm{coh}}(X)$ has a Serre
functor \cite{Bondal4}. This is presumably well-known.
\begin{propositions} Let $\Escr,\Fscr\in D^b_{\mathrm{coh}}(X)$. Then there
  are natural isomorphisms
\begin{equation}
\label{ref-6.1-54}
\Hom_{\Oscr_X}(\Escr,\Fscr)\cong \Hom_{\Oscr_X}(\Fscr,S\Escr)^\ast
\end{equation}
where $S\Escr=\Escr\otimes_{\Oscr_X} \omega_X[n]$.
\end{propositions}
\begin{proof}
We start with classical Serre duality \cite{ramis}:
\begin{itemize}
\item $H^n(X,\omega_X)=\CC$.
\item Let $\Fscr\in \coh(X)$. The Yoneda pairing
\[
H^i(X,\Fscr)\otimes \Ext^{n-i}_{\Oscr_X}(\Fscr,\omega_X)\r \CC
\]
is non-degenerate.
\end{itemize}
Now let $\Fscr\in D_{\mathrm{coh}}^b(\Oscr_X)$. {}{}From the pairing
\[
R\Gamma(X,\Fscr)\otimes_k R\Hom_{\Oscr_X}(\Fscr,\omega_X)\r
R\Gamma(X,\omega_X)\r \CC[-n]
\]
we obtain a map
\begin{equation}
\label{ref-6.2-55}
R\Gamma(X,\Fscr)\r  R\Hom_{\Oscr_X}(\Fscr,\omega_X[n])^\ast
\end{equation}
We claim that this is an isomorphism.
By induction over triangles we reduce to the case  $\Fscr\in \coh(X)$.
Then to show that \eqref{ref-6.2-55} is an isomorphism we have to show that
it is an isomorphism on cohomology, which is precisely classical Serre
duality.

If $\Escr\in D^-(X)$, $\Fscr\in D^+(X)$ then we
have the usual local-global isomorphism \cite{Verdier}:
\[
R\Hom_{\Oscr_X}(\Escr,\Fscr)=R\Gamma(X,\uRHom_{\Oscr_X}(\Escr,\Fscr))
\]
Now assume $\Escr,\Fscr\in D_{\mathrm{coh}}^b(X)$, $\Gscr\in D^+(X)$. We claim that
the following holds.
\begin{itemize}
\item[(a)]
$\uRHom_{\Oscr_X}(\Escr,\Fscr)\in D_{\mathrm{coh}}^b(X)$.
\item[(b)]
The natural map
\[
\uRHom_{\Oscr_X}(\Fscr,\Escr\Lotimes \Gscr)\r \uRHom_{\Oscr_X}(\uRHom_{\Oscr_X}(\Escr,\Fscr),\Gscr)
\]
is an isomorphism.
\end{itemize}
Since these statements are local we may assume that $\Escr$, $\Fscr$
are bounded free complexes. In that case (a) and (b) are obvious.

The proof of the proposition now follows from the following
computation:
\begin{align*}
R\Hom_{\Oscr_X}(\Escr,\Fscr)&=R\Gamma(X,\uRHom_{\Oscr_X}(\Escr,\Fscr))\\
&\cong R\Hom_{\Oscr_X}(\uRHom_{\Oscr_X}(\Escr,\Fscr),\omega_X[n])^\ast\\
&=R\Gamma(X,\uRHom_{\Oscr_X}(\uRHom_{\Oscr_X}(\Escr,\Fscr),\omega_X[n]))^\ast\\
&=R\Gamma(X,\uRHom_{\Oscr_X}(\Fscr,\omega_X[n]\otimes \Escr))^\ast\\
&=R\Hom_{\Oscr_X}(\Fscr,\omega_X[n]\otimes \Escr)^\ast \qed
\end{align*}
\def\qed{}\end{proof}
\subsection{Comparison of $\Ext$}
If $X$ is algebraic then it is
well-known and easy to prove that $D^b(\coh(X))$ and
$D^b_{\mathrm{coh}}(X)$  are equivalent. We don't know if the
corresponding result is true for
the complex analytic case.
For surfaces it is implied by the following proposition.
\begin{propositions}
\label{ref-6.2-56}
  Let $X$ be a smooth compact analytic  surface.
Then the Yoneda $\Ext$-groups in $\coh(X)$ coincide
with the $\Ext$-groups in the category of all $\Oscr_X$-modules.
\end{propositions}
\begin{proof}
Let us respectively write ${}^{I}\Ext$ and ${}^{II}\Ext$ for the
Yoneda $\Ext$ and the $\Ext$ in $\Mod(\Oscr_X)$.
Both $\Ext$'s are $\delta$-functors in their first and second argument
and they coincide in degree zero. Hence to
show that  ${}^{I}\Ext={}^{II}\Ext$  it is sufficient to show that
${}^{II}\Ext$ is elementwise effaceable in its first argument \cite[Lemma II.2.1.3]{Illusie}. That is if $i>0$,
$\Escr,\Fscr\in \coh(X)$ and $f\in \Ext^i(\Escr,\Fscr)$ then we have to show that there exists an epimorphism
$\Escr'\r\Escr$ in $ \coh(X)$ such that the image of $f$ under  the induced map ${}^{II}\Ext^i(\Escr,\Fscr)\r
{}^{II}\Ext^i(\Escr',\Fscr)$ is zero.

Let $\Escr,\Fscr\in\coh(X)$. We clearly have
\[
{}^{I}\Ext^1(\Escr,\Fscr)={}^{II}\Ext^1(\Escr,\Fscr)
\]
since the extension of two coherent sheaves is coherent. Since ${}^I\Ext^1$ is effaceable, so is  ${}^{II}\Ext^1$.

Furthermore we also have ${}^{II}\Ext^i(\Escr,\Fscr)=0$ for $i>2$.
This follows for example from \eqref{ref-6.1-54}.
Hence by \cite{BBD} we only have to
show that ${}^{II}\Ext^2$ is effaceable. To do this we use the
following sublemma:
\begin{sublemma}
Let $\Escr,\Fscr\in \coh(X)$. Choose $x\in X$ and let $m_x$ be the
corresponding maximal ideal in $\Oscr_X$. Then there exists $n$ such
that ${}^{II}\Ext^2(m_x^n\Escr,\Fscr)=0$.
\end{sublemma}
\begin{proof} By \eqref{ref-6.1-54} it suffices to show that for $n\gg 0$
  one has
  $\Hom(\Gscr,m_x^n\Escr)=0$ with $\Gscr=\Fscr\otimes
  \omega_X^{-1}$. Since $\Hom(\Gscr,m_x^n\Escr)$ is finite dimensional
  it is clearly sufficient to show that for $a\in \NN$ there exists
  $b>a$ such that $\Hom(\Gscr,m_x^b\Escr)\neq \Hom(\Gscr,m_x^a\Escr)$.

So pick a non-zero $f:\Gscr\r m_x^a\Escr$. Then there will exist $b$ such that
$\im f_x\not\subset m_x^b\Escr_x$ (look at stalks). Hence $f\not\in
\Hom(\Gscr,m_x^b\Escr)$. This finishes the proof.
\end{proof}

To complete the proof that ${}^{II}\Ext^2$ is effaceable we pick $x\neq
y$ in $X$ and we choose $n$ such that ${}^{II}\Ext^2(m_x^n\Escr,\Fscr)=
{}^{II}\Ext^2(m_y^n\Escr,\Fscr)=0$. Since the canonical map
$m_x^n\Escr\oplus m_y^n\Escr\r \Escr$ is surjective, we are done.
\end{proof}
\begin{corollarys}
\label{6.3aa}
Let $X$ be as above. Then the canonical functor $F:D^b(\coh(X))\r D^b_{\mathrm{coh}}(X)$ is
  an equivalence.
\end{corollarys}
\begin{proof} By induction over triangles and the above proposition we
  see that $F$ is fully faithful. That it is essentially surjective
  also follows by induction over triangles.
\end{proof}

\subsection{The derived category of an exact category}
Assume that $\Escr$ is an exact category \cite{Quillen}. In \cite{Neeman4}
Neeman defines the derived category $D(\Escr)$ of $\Escr$.  By
definition $D(\Escr)=K(\Escr)/K(\Escr)^{eac}$ where as usual
$K(\Escr)$ is the homotopy category of $\Escr$ and  $K(\Escr)^{eac}$
is the epaisse envelope of the category
$K(\Escr)^{ac}$ of acyclic 
complexes in $K(\Escr)$.  By definition a complex 
\[
\cdots\r X^n\r X^{n+1}\r X^{n+2}\r\cdots
\]
is \emph{acyclic} if each map $X^n\r X^{n+1}$ decomposes in $\Escr$ as a
composition of an admissible epimorphism with an admissible
monomorphism: $X^n\r
D^n\r X^{n+1}$ such that $D^n\r X^{n+1}\r D^{n+1}$ is exact.  Since by
\cite[Lemma 1.1]{Neeman4} $K(\Escr)^{ac}$  is triangulated it follows from 
lemma \ref{epaisseenvelope} that every object in $K(\Escr)^{eac}$ is a direct
summand of an object in $K(\Escr)^{ac}$. Furthermore if $\Escr$ is
Karoubian then by \cite[Lemma 1.2]{Neeman4}
$K(\Escr)^{eac}=K(\Escr)^{ac}$.
\subsection{Torsion pairs in abelian categories}
Assume that $\Cscr$ is an abelian category and let $(\Tscr,\Fscr)$ be
a {\em torsion pair} in $\Cscr$, i.e. $\Tscr$ and $\Fscr$ are full
subcategories in $\Cscr$
such that
$\Hom(\Tscr,\Fscr)=0$ and every  object
$C\in\Cscr$ fits in a unique exact sequence
\begin{equation}
\label{tcf}
0\r T\r C\r F\r 0
\end{equation}
with $T\in \Tscr$ and $F\in \Fscr$.  It follows that $\Tscr$ and
$\Fscr$ are respectively closed under quotients and subobjects.

The assignments $C\mapsto T$ and $C\mapsto F$ in the exact sequence \eqref{tcf}
yield functors $\tau : \Cscr\r \Tscr$ and $\phi : \Cscr \r \Fscr$ 
which are respectively the right and left adjoint to the inclusions $\Tscr\r
\Cscr$, $\Fscr\r \Cscr$.  

%  The inclusions $\Fscr\subset
% \Cscr$, $\Tscr\subset \Cscr$ are respectively left and right exact.
% The adjoints
%  $\phi$ and $\tau$ are respectively right and left exact.

It is easy to see that $\Tscr$ and $\Fscr$ possess  kernels and
 cokernels. We have formulas
\begin{equation}
\label{formulas}
\begin{split}
\ker_{\Fscr}&=\ker_{\Cscr}\\
\coker_{\Fscr}&=\phi\circ \coker_\Cscr
\end{split}
\end{equation}
and dual formulas for $\Tscr$.

Following \cite{HRS} we say that
$(\Tscr,\Fscr)$ is {\em tilting} if every object in $\Cscr$ is a subobject
of an object in $\Tscr$. Similarly $(\Tscr,\Fscr)$ is {\em cotilting} if
every object in $\Cscr$ is a quotient of an object in $\Fscr$.

The torsion pair $(\Tscr,\Fscr)$ defines a $t$-structure on
$D^b(\Cscr)$ by
\begin{align*}
{}^pD^b(\Cscr)^{\le 0}&=\{C\in D^b(\Cscr)^{\le 1}\mid H^1(C)\in\Tscr\}\\
{}^pD^b(\Cscr)^{\ge 0}&=\{C\in D^b(\Cscr)^{\ge 0}\mid H^0(C)\in\Fscr\}
\end{align*}
By definition the tilting ${}^p\Cscr$ of $\Cscr$ with respect to
$(\Tscr,\Fscr)$ is the heart of this $t$-structure. It is easy to see
that $(\Fscr,\Tscr[-1])$ is a torsion pair in ${}^p\Cscr$. Furthermore
according to \cite[Prop. I.3.2]{HRS} $(\Tscr,\Fscr)$ is tilting if and
only if $(\Fscr,\Tscr[-1])$ is cotilting and vice versa.

\smallskip

Let $\Escr$ be either $\Tscr$ or $\Fscr$.
The exact structure on $\Cscr$ induces
an exact structure on $\Escr$. This is 
intrinsically determined in the following way:
 a morphism $f:A\r B$ in $\Escr$ is \emph{strict} if the canonical morphism $\coker 
\ker f\r \ker \coker f$ is an isomorphism.  A diagram
\[
0\r A\xrightarrow{f} B\xrightarrow{g} C\r 0
\]
is an admissible exact sequence if $f$ is a strict monomorphism, $g$
is a strict epimorphism and $\coker f=g$, $\ker g=f$.

The following statements are obvious.
\begin{lemmas} 
\label{obvious}
\begin{enumerate}
\item
A complex over $\Escr$ is acyclic if and only if it is
  acyclic in~$\Cscr$.
\item $K(\Escr)^{eac}=K(\Escr)^{ac}$. 
\item A map between complexes over $\Escr$ is an isomorphism in
  $D(\Escr)$ if and only if it is a quasi-isomorphism over $\Cscr$. 
\end{enumerate}
\end{lemmas} 
\begin{lemmas} 
\label{cotilting}
\cite[Ex.\ 1.3.23(iii)]{BBD} Assume that $(\Tscr,\Fscr)$ is cotilting. Then the canonical
  map $D(\Fscr)\r D(\Cscr)$ is an equivalence.
\end{lemmas}
\begin{proof}  Since $(\Tscr,\Fscr)$ is cotilting and $\Fscr$ is
  closed under subobjects, every object in $\Cscr$ has
  a resolution of length two by objects in $\Fscr$. Therefore by the
  (dual version) of \cite[Lemma I.4.6]{RD}  it follows that if $X$ is a 
  complex over $\Cscr$  there exists a quasi-isomorphism $F\r X$
  with $F$ a complex over $\Fscr$.

We find for $F_1,F_2$ complexes over $\Fscr$
\begin{align*}
\Hom_{D(\Cscr)}(F_1,F_2)&=\dirlim_{X\xrightarrow[\mathrm{qi}]{} F_1}
\Hom_{K(\Cscr)}(X,F_2)\\
&=\dirlim_{\begin{smallmatrix} F'_1\xrightarrow[\mathrm{qi}]{} F_1\\ F'_1\in K(\Fscr)
    \end{smallmatrix}} 
\Hom_{K(\Cscr)}(F'_1,F_2)\\
&=\Hom_{D(\Fscr)}(F_1,F_2)
\end{align*}
The last equality follows from lemma \ref{obvious}(3).
\end{proof}
This result was also proved by Schneiders in the (equivalent) setting
of quasi-abelian categories. This is explained in Appendix \ref{quasiabelian}. 

The following result is proved in \cite{HRS} under some
additional (unnecessary) conditions.
\begin{propositions}
\label{ref-6.3.1-57}
Assume that $(\Tscr,\Fscr)$ is cotilting. Then
  $D({}^p\Cscr)=D(\Cscr)$.
\end{propositions}
\begin{proof}
According to lemma \ref{cotilting}        we have $D(\Cscr)=D(\Fscr)$. Since
$(\Fscr,\Tscr[-1])$ is tilting, we can invoke the dual result for
${}^p\Cscr$ which is $D({}^p\Cscr)=D(\Fscr)$.  Since the
exact structure on $\Fscr$ is intrinsic, the induced exact structures
on $\Fscr$ from the inclusions $\Fscr\r \Cscr$ and $\Fscr\subset
{}^p\Cscr$ are the same
and this finishes the proof.
\end{proof}

\begin{remarks}  Lemma \ref{cotilting} and
  Propositions \ref{ref-6.3.1-57} are also valid for 
  $D^b$ (the equivalences preserve boundedness).
\end{remarks} 

\subsection{Tilting in noetherian abelian categories}
\begin{lemmas}
%\label{ref-6.3.1-57}
\label{noeth}
Let $(\Tscr,\Fscr)$ be a torsion pair in $\Cscr$.
Then $\Cscr$ is noetherian if and only if the
following conditions hold:
\begin{enumerate}
\item[N1.]
Every chain of subobjects of $F$: $F_0\subset F_1\subset F_2\subset \cdots$ for
$F_i\in \Fscr$, $F\in\Fscr$ becomes stationary.
\item[N2.]
Every chain of epimorphisms $T_0\r T_1\r T_2\r \cdots$ for
$T_i\in\Tscr$ becomes
stationary.
\end{enumerate}
\end{lemmas}
\begin{proof}
Let us show that N1, N2 imply $\Cscr$ noetherian. Let $C_0\subset C_1\subset
 \cdots$ be an ascending chain of subobjects of $C\in \Cscr$.
The sequence $F_i={\rm Im}(\phi (C_i)\to \phi (C))$ becomes stationary by
N1. Denote by $F\subset \phi (C)$ the limiting subobject of the sequence.
We may relace $C$ by the fibred product $C'=C\times _{\phi (C)}F={\rm
ker}(C\oplus F\r \phi (C))$. Indeed, the natural morphisms $C_i\r C'$ are
monic, because so are the composites $C_i\r C'\r C$.

By construction of $C'$,
$\phi(C')=F$ and the maps $\phi (C_i)\r \phi (C')$ are epic for $i\gg 0$. If
$R_i=C'/C_i$, then we have a complex
$\phi (C_i)\to \phi (C')\to \phi (R_i)$.
As $\phi $ is a left adjoint it takes epi to epi, so both morphisms in this complex are epic. It follows
that $\phi (R_i)=0$, i.e. $R_i\in \Tscr$ for $i\gg 0$.
Therefore, the chain of epimorphisms
$R_0\to R_1\to \dots $ becomes stationary by N2. This proves that
the primary chain of $C_i$'s becomes stationary. The converse statement is
obvious.
\end{proof}

By \eqref{formulas}
morphisms in $\Tscr$ are epimorphisms iff they are epimorphisms
in $\Cscr$
and morphisms in $\Fscr$ are monomorphisms iff they are monomorphisms
in $\Cscr$.
So N1 and N2 are intrinsic in $\Tscr,\Fscr$.

We will use the following criterion for ${}^p\Cscr$ to be noetherian.
\begin{lemmas}
\label{ref-6.3.2-58}
 Assume that $\Cscr$ is noetherian and $(\Tscr,\Fscr )$ a torsion pair in
$\Cscr$. Then ${}^p\Cscr$ is
  noetherian if and only if the following is true: every ascending
  chain $F_0\subset F_1\subset \cdots $ with $F_i\in \Fscr$ and
  $\coker(F_0\r F_i)\in\Tscr$ for all $i$, is stationary.
\end{lemmas}
\begin{proof}
If there is an ascending chain as in the statement of the lemma which
is not stationary then it is easy to see that we have an ascending chain of subobjects of $F_0$ in ${}^p\Cscr$
\[
F_1/F_0[-1]\subset F_2/F_0[-1]\subset F_3/F_0[-1]\subset \cdots.
\]
 Hence ${}^p\Cscr$ is not noetherian. So we will now concentrate on the
converse direction.

%It is easy to prove that $\Cscr$ will be noetherian if and only if the
%following conditions hold:
%\begin{enumerate}
%\item
%Every chain of subobjects of $F$: $F_0\subset F_1\subset F_2\subset \cdots$ for
%$F_i\in \Fscr$, $F\in\Fscr$ becomes stationary.
%\item
%Every chain of epimorphisms $T_0\r T_1\r T_2\r \cdots$ for
%$T_i\in\Tscr$ becomes
%stationary.
%\end{enumerate}
%Since morphisms in $\Tscr$ are epimorphisms iff they are epimorphisms
%in $\Cscr$
%and morphisms in $\Fscr$ are monomorphisms iff they are monomorphisms
%in $\Cscr$,
%1. and 2. are intrinsic in $\Tscr,\Fscr$.

By lemma \ref{noeth}, to check that ${}^p\Cscr$ is noetherian we have to verify
N1, N2 with
$\Tscr$ and $\Fscr$ exchanged. To this end we have to know the nature
of monomorphisms in $\Tscr$ and epimorphisms in $\Fscr$. {}From
\eqref{formulas} we obtain:
\begin{itemize}
\item Monomorphisms in $\Tscr$ are the maps whose kernel in $\Cscr$
  is in $\Fscr$.
\item Epimorphisms in $\Fscr$ are the maps whose cokernel in $\Cscr$
  is in $\Tscr$.
\end{itemize}
Let us now check that N2 holds if we replace $\Tscr$ by $\Fscr$.
%We assume that $\Cscr$ is noetherian.
Thus we have a chain of maps in $\Fscr$
\begin{equation}
\label{ref-6.3-59}
F_0\r F_1\r F_2 \r \cdots
\end{equation}
whose cokernel is in $\Tscr$.
Using the fact that $\Cscr$ is
noetherian we see that the kernel $K_{ij}=\ker(F_i\r F_j)$ will become
stationary for $j\gg 0$.  Let $K_i=K_{ij}$ for $j\gg 0$. Then the maps
$F_i/K_{i}\r F_{i+1}/K_{i+1}$ are injective. Using the fact that
$F_i/K_i$ injects in $F_j$ for $j\gg 0$ we see that
$F_i/K_i\in\Fscr$.  Furthermore $(F_{i+1}/K_{i+1})/(F_i/K_i)$ is a
quotient of $\coker (F_i\r F_{i+1})$ so it lies in $\Tscr$.

It follows that the condition given in the statement of the lemma
holds for the sequence $(F_i/K_i)_i$, i.e. this sequence  will become stationary. Hence by left
shifting if necessary we may assume that $F_i/K_i\r  F_{i+1}/K_{i+1}$
is an isomorphism for all $i\ge 0$. {}From the snake lemma we then deduce that
$\coker(K_i\r K_{i+1})$ is isomorphic to $\coker(F_i\r F_{i+1})$ and
hence is in $\Tscr$. This implies that for $j\gg 0$, $K_j=\coker(K_0\r
K_j)\in\Tscr$. Since also $K_j\in\Fscr$ this implies $K_j=0$ for $j\gg
0$.
Truncating the beginning of the sequence by sufficiently big $j$ we obtain a sequence which
satisfies the conditions in the statement of the lemma.
This implies that $F_j\r F_{j+1}$ is an isomorphism
for $j\gg 0$.

Let us now assume N2
and check that N1 holds if we replace $\Fscr$ by $\Tscr$.  Thus
we have a chain of maps in $\Tscr$
\[
T_0\r T_1\r T_2 \r \cdots \r T
\]
whose kernel is in $\Fscr$. Since $\Cscr$ is noetherian, the images of
the maps $T_i\r T$ will become stationary. Since these images are in
$\Tscr$ we may without loss of generality assume that the maps $T_i\r
T$ are surjective. Put $F_i=\ker(T_i\r T)$.
%By the snake lemma we obtain that the maps
Then $\coker(F_i\r F_{i+1})$, being isomorphic
to $\coker(T_i\r T_{i+1})$, is in
$\Tscr$. Hence
the chain $(F_i)_i$ is like that in \eqref{ref-6.3-59}, hence it
becomes stationary.
%But using the snake lemma again t
This implies that the chain $(T_i)_i$ also becomes
stationary.
\end{proof}
\begin{remarks} If $\Tscr \subset \Cscr$ is the subcategory of torsion sheaves
in the category of coherent sheaves on an analitic or algebraic variety
(the case of our interest in the next subsection), then $\Tscr $ has a
property to be closed under subobjects in $\Cscr$. Under this additional
condition the proof of the lemma can be simplified in two places:
$\coker (K_i\r K_{i+1})$ are torsion being subobjects of $\coker (F_i\r
F_{i+1})$ and N1 with $\Fscr $ replaced by $\Tscr $ is automatically satisfied
as the kernels of $T_i\r T_{i+1}$ and $T_i\r T$ are trivial.
\end{remarks}

\subsection{Non-saturation for analytic surfaces}
We can now prove the following result:
\begin{theorems}
Let $X$ be a smooth compact analytic surface with no curves. Then $D^b(\coh(X))$ is not
saturated.
\end{theorems}
By Corollary \ref{6.3aa} the result also holds for $D^b_{\coh}(X)$.
\begin{proof}
\begin{step}
\label{ref-1-60}
Let $\Tscr \subset\coh(X)$ be the full subcategory of objects in $\coh(X)$ whose
  support is strictly smaller than $X$. Since $X$ contains no curves
  and is compact, this support must be a finite set of points.
Let $\Fscr$ be the full subcategory of objects $F$ in $\coh(X)$ such that
  $\Hom(\Tscr,F)=0$.  It is clear that $(\Tscr,\Fscr)$ is a torsion
  pair.
\end{step}
\begin{step}
\label{ref-2-61}
$\Tscr$ is closed under essential extensions.
  To prove this let $T\in\Tscr$ and let $T\subset T'$ be an essential
  extension. Let $\{x_1,\ldots,x_n\}\in X$ be the support of $T$.
By the Artin-Rees property of the stalks of $\Oscr_{X,x_i}$ there
  exists $t\ge 0$ such that $m_{x_i}^tT'_{x_i}\cap T_{x_i}=0$ for all
  $i$. Thus $m_{x_1}^t\cdots m_{x_n}^tT'\cap T=0$ and since we are in an
  essential extension it follows $m_{x_1}^t\cdots
  m_{x_n}^tT'=0$. Hence $T'\in\Tscr$.
\end{step}
\begin{step}
\label{ref-3-62}
Let $E\in \coh(X)$. Then $E$  is a quotient
  of an object in $\Fscr$. In \cite{schuster} Schuster proved the more general result that every coherent sheaf on a complex surface is a quotient of a vector bundle. We give a simple proof of the weaker statement that we need.

We write $E$ as an extension
\[
0\r T\r E\r F\r 0
\]
where $T\in\Tscr$ is torsion and $F\in\Fscr$.
Take the maximal $E'\subset E$, such that $E'\cap T=0$. As $T\subset E/E'$ is an essential extension, then by the previous step $E/E'\in\Tscr$. $E'\in\Fscr$ by the choice of $T$. We now obtain and exact
sequence
\[
0\r E'\r E\r T'\r 0
\]
with $T'\in \Tscr$. It is easy to see that every object in $\Tscr$ is
a quotient of a free $\Oscr_X$-module. So write $T'$ as a quotient of $F'\in
\Fscr$ and let $E''$ be the corresponding pullback of $E$. Then
$E''$ is an extension of $E'$ and $F'$ and hence $E''\in \Fscr$. Thus
we have written $E$ as a quotient of $E''\in\Fscr$.
\end{step}
\begin{step}
\label{ref-4-63}
By the previous step $(\Tscr,\Fscr)$ is cotilting. Hence by
lemma \ref{ref-6.3.1-57} $D^b(\coh(X))=D^b({}^p\!\coh(X))$.
\end{step}
\begin{step}
\label{ref-5-64}
Now we claim that ${}^p\!\coh(X)$ is noetherian. By lemma
  \ref{ref-6.3.2-58} we need to show that every ascending chain
\[
F_0\subset F_1\subset F_2\subset\cdots
\]
with $F_i\in\Fscr$, $F_i/F_0\in\Tscr$
 becomes stationary.

This is satisfied in our case because we must have $F_n\subset
F_0^{\ast\ast}$ and $F_0^{\ast\ast}/F_0$ has finite length.
\end{step}
\begin{step}
\label{ref-6-65}
Note that
  ${}^p\!\coh(X)$ is self-dual under $R\uHom(-,\Oscr_X)$. Hence it is both
  noetherian and Artinian. Thus ${}^p\!\coh(X)$ has finite length.
\end{step}
\begin{step}
\label{ref-7-66}
Assume that $\coh(X)$ is saturated. By Step \ref{ref-4-63}
  ${}^p\!\coh(X)$ will also be saturated. Since this is a finite length
  category it follows from lemma \ref{ref-2.1-22} that it has to be of
the form $\mod(\Lambda)$ for a
  finite dimensional algebra $\Lambda$.

%For the surfaces we consider $\omega_X$ is trivial and hence the Serre functor
%is a pure shift by two places.
By proposition  \ref{ref-6.1-54},  $S[-2]$ is a functor which preserves $\coh(X)$, $\Tscr$ and $\Fscr$. Hence it
preserves ${}^p\!\coh(X)$ (regarded as a subcategory in $D^b(\coh(X)$). For
a finite dimensional algebra the Serre functor takes projectives into
injectives. Therefore, its shift by $-2$ cannot preserve the category
$\mod(\Lambda)$. We have obtained a contradiction.

%Since for a finite dimensional algebra the
%  Serre functor can never be a non-trivial shift
%  (for example a projective goes to an injective) we have obtained a
%contradiction. \qed
\end{step}
\def\qed{}\end{proof}
\begin{remarks} It seems likely that this counter example is only the
  tip of the iceberg and
  that in fact a compact analytic manifold is saturated if and only if
  it is an algebraic space. This would mean that saturatedness would be
  a criterion for a triangulated category to be of algebraic
  nature.
\end{remarks}
Among the surfaces to which the theorem is applicable are K3, 2-dimensional
tori and surfaces of type VII in the Kodaira classification \cite{Nakamura}.

\appendix
\section{An alternative proof in the commutative case}
Theorem \ref{ref-1.1-1}, as stated follows from the non-commutative
result Theorem \ref{ref-5.2.12-53}. However in the commutative case it
is possible to give a straightforward proof of a more general
result.
\begin{theorem}
  Assume that $X$ is a projective variety over a field $k$. Let
  $\Dscr$ be the triangulated category of perfect complexes on $X$. Then every
  contravariant cohomological functor of finite type on $\Dscr$ is
  representable by a bounded complex with coherent homology.
\end{theorem}
\begin{proof} According to \cite[lemma 2.13]{keller3} $H$ is
  represented by an object  $E$ in $D(\Qch(X))$. We have to show that
  this object is in $D^b(\coh(X))$. To prove this we repeat the
  argument of \cite{Bondal4}.

Choose an embedding $\pi:X\r \PP^n$ and consider the functor
$H'=H\circ R\pi^\ast:D^b(\coh(\PP^n))\r \Vect(k)$. According to
Beilinson's result \cite{Beilinson} as it was reformulated in \cite{Baer,Bondal2}
there is an equivalence $\theta:D^b(\mod(\Lambda))\r D^b(\coh(\PP^n))$
where $\Lambda$ a finite dimensional algebra of finite global
dimension. Put $H''=H'\circ \theta$.  Invoking \cite[lemma 2.13]{keller3} again we
see that $H''$ is representable by an object $G$ in
$D(\Lambda)$.  Since $H''$ is still of finite type it follows that $\sum_n \dim H^{\prime\prime n}(\Lambda)=\sum_n
\dim \Hom(\Lambda[n],G)<\infty$. Thus $G\in D^b(\mod(\Lambda))$.
 This
implies that  $H'$ is represented by $F=\theta^{-1}(G)\in D^b(\coh(X))$.

Thus if $A\in D^b(\coh(\PP^n))$  we have
\begin{align*}
\Hom_{\PP^n}(A,R\pi_\ast E)&=\Hom_X(R\pi^\ast A,E)\\
&=H(R\pi^\ast A)\\
&=H'(A)\\
&=\Hom_{\PP^n}(A,F)
\end{align*}
Putting $A=F$ we obtain a map $\mu:F\r \pi_\ast E$ which becomes an
isomorphism if we apply $\Hom_{\PP^n}(A,-)$ for $A\in D^b(\coh(\PP^n))$.  In other
words the cone of $\mu$ is right orthogonal to
$D^b(\coh(\PP^n)))$. By taking $A=\Oscr(n)_n$ it follows easily that
the cone of $\mu$ is zero and hence $\mu$ is an isomorphism.  Thus
$\pi_\ast E\in D^b(\coh(\PP^n))$. This implies $E\in D^b(\coh(X))$
\end{proof}
\section{Quasi-abelian categories}
\label{quasiabelian}
In this appendix we discuss quasi-abelian categories. Let $\Escr$ be
an addititive category with kernels and cokernels.  A morphism $f:A\r
B$ is said to be \emph{strict} if the canonical map $\coker \ker f\r
\ker \coker f$ is an isomorphism.

We say $\Escr$ is quasi-abelian if $\Escr$ satisfies the property that
the pullback of any strict epi is strict epi and the pushout of any
strict mono is strict mono. Quasi-abelian category
appear frequently in the literature, often under different names.
They are called ``preabelian'' in \cite{Jurchescu}, ``semiabelian'' in
\cite{Raikov} and quasi-abelian in \cite{Schneiders,Cruciani}. It can
also be seen that  quasi-abelian categories are additive categories which
are regular and coregular \cite{BGO}.
In this appendix we show that
the notion of a quasi-abelian category 
is the same as that of a (co)tilting torsion theory.

\medskip

Let $\Escr$ be quasi-abelian. $\Escr$ carries an intrinsic exact
structure with the admissible mono- and epimorphism being
respectively the strict mono- and epimorphisms \cite[\S 1.1.4]{Schneiders}.

In \cite[\S 1.2.3]{Schneiders} it is shown that $\Escr$ has two canonical embeddings
into abelian categories $\Lscr\Hscr(\Escr)$ and $\Rscr\Hscr(\Escr)$
preserving and reflecting exactness. Furthermore $\Escr$ is stable
under extensions in these embeddings.

\begin{proposition} \label{embedding} \cite[Prop. 1.2.35]{Schneiders} The embedding $\Escr\subset
\Lscr\Hscr(\Escr)$ is characterized by the following properties: 
$\Escr\subset \Lscr\Hscr(\Escr)$ is a fully faithful embedding of $\Escr$ into
an abelian category, $\Escr$ is closed under subobjects in
$\Lscr\Hscr(\Escr)$  and  every object in $\Lscr\Hscr(\Escr)$ is a quotient of an object in
$\Escr$. 
\end{proposition}
The following result is  \cite[Prop. 1.2.31]{Schneiders}.
\begin{proposition} \label{derivedeq}
The inclusion $\Escr\subset \Lscr\Hscr(\Escr)$ extends to an
equivalence of derived categories $D(\Escr)\cong D(\Lscr\Hscr(\Escr))$.
\end{proposition}
The following result shows that the notion of a quasi-abelian category
is the same as that of (co)tilting torsion theory.
\begin{proposition} 
\label{equivalence} Let $\Escr$ be an additive category. The following
  are equivalent.
\begin{enumerate} 
\item $\Escr$ is quasi-abelian.
\item There exists a cotilting torsion pair $(\Tscr,\Fscr)$ in an
  abelian category $\Cscr$ with $\Escr\cong \Fscr$.
\item There exists a tilting torsion pair $(\Tscr',\Fscr')$ in an abelian
  category $\Cscr'$ with $\Escr\cong \Tscr'$.
\end{enumerate}
In the situation of (2) we have $\Cscr\cong \Lscr\Hscr(\Fscr)$ and in
the situation of (3) we have $\Cscr'\cong \Rscr\Hscr(\Fscr)$.
\end{proposition}
\begin{proof}
That $\Cscr\cong \Lscr\Hscr(\Escr)$ and $\Cscr'\cong
\Rscr\Hscr(\Escr)$ follows directly from Proposition \ref{embedding}
(and its dual version).

To prove the stated equivalence we note that by symmetry we only need to prove the equivalence of (1) and (2).
\begin{itemize}
\item[(2)$\Rightarrow$(1)]  Since $\Fscr$ is exact, pullbacks of
  admissible epimorphisms are admissible epimorphisms. Since the
  admissible epimorphisms are precisely the strict epimorphisms this
  shows that pullbacks of strict epimorphisms are strict
  epimorphisms. The corresponding result for strict monomorphisms is
  proved in the same way.

% We have already noted that $\Fscr$ has
%   kernels and cokernels. Assume first that $f:F\r G$ is a
%   strict epimorphism in $\Fscr$ and $H\r G$ is an arbitrary morphism
%   in $\Fscr$. Let $f':P\r H$ be
%   the pullback of $f$ in $\Cscr$.  Since $P\subset H\oplus F$ it
%   follows that $P\in \Fscr$ and hence $f'$
% is also the
%   pullback of $f$ in $\Fscr$.

% We need to show that $f'$ is a strict
%   epimorphism. {}From the formulas \eqref{formulas} it is easy to see that the
%   strict epimorphisms in $\Fscr$ correspond precisely to the epimorphisms in
%   $\Cscr$. Hence the result we need to show follows from the
%   corresponding result for abelian categories.

% Now assume that $r:R\r S$ is a strict monomorphism and let $R\r T$ be an
% arbitrary morphism. Let $r':T\r Q$ be the pushout of $r$ in $\Cscr$.  It
% is clear that $r'$ is a monomorphism  in $\Cscr$. 

% Again using the formulas \eqref{formulas} it is easy to see that
% a morphism in $\Fscr$  is a strict monomorphism in $\Fscr$ if and only
% if it is monomorphism in $\Cscr$ and if its cokernel in $\Cscr$ lies  in $\Fscr$.

% Now we have $\coker_\Cscr r'\cong \coker_\Cscr r\in \Fscr$ since $r$ is strict.
% Thus $Q$ is an extension of objects in $\Fscr$ and hence it is itself
% in $\Fscr$. Thus $r'$ is also the pushout of $r$ in $\Fscr$. Since
% $\coker_\Cscr r'\in \Fscr$ it also follows that $r'$ is strict.

\item[(1)$\Rightarrow$(2)]
Put $\Fscr=\Escr$ and $\Cscr=\Lscr\Hscr(\Escr)$. Let $\Tscr$ be the
full subcategory of $\Cscr$ consisting of objects $\coker_\Cscr f$
where the morphism $f$ is  an epimorphism in $\Fscr$. 

We claim that $(\Tscr,\Fscr)$ is a cotilting torsion pair in $\Cscr$. 
If  $T=\coker_\Cscr f \in \Tscr$ and  $F\in \Fscr$ then from the fact that $f$ is an
epimorphism in $\Fscr$ we immediately obtain $\Hom(T,F)=0$.

Now let $C$ be an arbitrary object in $\Cscr$. According to
Proposition \ref{embedding} there exist a short exact sequence in $\Cscr$
\begin{equation}
F \xrightarrow{f} F' \r C\r 0
\end{equation}
with $F,F'\in\Cscr$.
In particular if $(\Tscr,\Fscr)$ is a torsion theory then it will
certainly be cotilting. 

We will now show that $C$ is an extension of the form \eqref{tcf}.
We have the following
commutative diagram
\begin{equation}
\label{comdiag1}
\begin{CD}
F @>f>> F' @>g'>> \coker_\Fscr f \\
 @V\alpha VV    @| @| \\
 \ker_\Fscr g' @>>f'> F' @>>g'> \coker_\Fscr f 
\end{CD}
\end{equation}
It is  easily checked that $\coker_\Fscr f$ satifies the universal
property for being a cokernel of $f'$. Thus $\coker_\Fscr f'=\coker_\Fscr f$ and
hence $g'$ a strict epimorphism. 

Hence we obtain in particular the following: \emph{a cokernel of an
  arbitrary morphism in $\Fscr$ is a strict epimorphism}. Dually we
  also obtain: 
  \emph{a kernel of an arbitrary morphism in $\Fscr$ is a strict
  monomorphism}. Thus in particular  $f'$ is a strict monomorphism. It also follows that
  the lower sequence in \eqref{comdiag1} is an admissible exact
  sequence.

 We claim that $\alpha$ is an epimorphism in $\Fscr$. To show this assume that
 there is a morphism $
\beta:\ker_\Fscr g'\r Z$ in $\Fscr$ whose composition with
 $\alpha$ is zero. We have to prove $\beta=0$.

We extend the commutative diagram  \eqref{comdiag1}
 as follows:
\begin{equation}
\label{comdiag2}
\begin{CD}
F @>f>> F' @>g'>> \coker_\Fscr f \\
 @V\alpha VV    @| @| \\
 \ker_\Fscr g' @>>f'> F' @>>g'> \coker_\Fscr f \\
@V\beta VV @VV\gamma V\\
Z @>f''>> Z'
\end{CD}
\end{equation}
where the lower square is a pushout in $\Fscr$. We now have $f''\circ
\beta\circ \alpha=0$ and hence $\gamma\circ f=0$. Thus $\gamma=\phi\circ
g'$ for some morphism $\phi:\coker_\Fscr f\r Z'$.

We deduce $f''\circ \beta=\gamma\circ f'=\phi\circ
g'\circ f'=0$.  Since we had assumed that
$\Fscr=\Escr$ is quasi-abelian we know that $f''$ is a strict
monomorphims and in particular a monomorphism. Thus it follows that $\beta=0$ and
hence $\alpha$ is an epimorphism.

Furthermore  by looking at the
decomposition 
\[
F\xrightarrow{\alpha} \ker_\Fscr g' \xrightarrow{f'} F'
\]
of $f$ in $\Cscr$ we find that
$C=\coker_{\Cscr} f$ is an extension of $\coker_\Cscr f'$ by
$\coker_\Cscr \alpha$. {}From the fact that $\alpha$ is an epimorphism in $\Fscr$ we obtain that
$\coker_{\Cscr} \alpha$ is in $\Tscr$. 
Now
since the lower sequence in \eqref{comdiag1} is an admissible exact
sequence and the embedding of $\Fscr\subset\Cscr$ preserves
exactness, we have $\coker_\Cscr f'=\coker_{\Fscr} f'\in \Fscr$.
This finishes the proof of (1)$\Rightarrow$(2).\qed
\end{itemize}
\def\qed{}\end{proof}
\begin{corollary} If $(\Tscr,\Fscr)$ is a cotilting torsion theory in
  an abelian category $\Cscr$ then $\Cscr\cong\Lscr\Hscr(\Fscr)$ and
  ${}^p\Cscr\cong\Lscr\Rscr(\Fscr)$.
\end{corollary}
\begin{proof}  By Proposition \ref{equivalence} we have
  $\Cscr=\Lscr\Hscr(\Fscr)$. Now $(\Fscr,\Tscr[-1])$ is a tilting torsion pair
  in ${}^p\Cscr$ and hence, again by Proposition \ref{equivalence}, we have
  $\Cscr=\Rscr\Hscr (\Fscr)$.
\end{proof}
Hence we find that Lemma \ref{cotilting} follows from Proposition
\ref{derivedeq}.

%\bibliography{mybibs}
%\bibliographystyle{amsabbrv}
\ifx\undefined\bysame
\newcommand{\bysame}{\leavevmode\hbox to3em{\hrulefill}\,}
\fi

\end{document}